\documentclass[12pt]{article}

\usepackage{amsmath, amssymb, booktabs, color, algorithm, algorithmic}
\usepackage{algorithm}

\usepackage{listings,color,url}

\definecolor{string}{rgb}{0,0.5,0.25}
\definecolor{comment}{rgb}{0,0.5,0}
\newtheorem{thm}{Theorem}[section]
\newtheorem{prop}{Proposition}[section]
\newtheorem{lem}{Lemma}[section]

\newtheorem{rem}{Remark}[section]

\newtheorem{prob}{Problem}[section]

\newtheorem{assum}{Assumption}[section]

\newcommand{\argmax}{\operatornamewithlimits{argmax}}
\newcommand{\argmin}{\operatornamewithlimits{argmin}}

\numberwithin{equation}{section}

\begin{document}
\makeatletter

\begin{center}
%%%%Performance Optimization of a Dynamic Channel Bonding Strategy in Cognitive Radio Networks
\large{\bf Two Stochastic Optimization Algorithms for Convex Optimization with Fixed Point Constraints}\\
\small{This work was supported by the Japan Society for the Promotion of Science through a Grant-in-Aid for
Scientific Research (C) (15K04763).}
\end{center}\vspace{2mm}

\begin{center}
\textsc{Hideaki Iiduka}\\
Department of Computer Science, 
Meiji University,
1-1-1 Higashimita, Tama-ku, Kawasaki-shi, Kanagawa 214-8571 Japan\\ 
(iiduka@cs.meiji.ac.jp)
\end{center}

\vspace{2mm}

\footnotesize{
\noindent\begin{minipage}{14cm}
{\bf Abstract:}
Two optimization algorithms are proposed for solving a stochastic programming problem for which the objective function is given in the form of the expectation of convex functions and the constraint set is defined by the intersection of fixed point sets of nonexpansive mappings in a real Hilbert space. This setting of fixed point constraints enables consideration of the case in which the projection onto each of the constraint sets cannot be computed efficiently. Both algorithms use a convex function and a nonexpansive mapping determined by a certain probabilistic process at each iteration. One algorithm blends a stochastic gradient method with the Halpern fixed point algorithm. The other is based on a stochastic proximal point algorithm and the Halpern fixed point algorithm; it can be applied to nonsmooth convex optimization. Convergence analysis showed that, under certain assumptions, any weak sequential cluster point of the sequence generated by either algorithm almost surely belongs to the solution set of the problem. Convergence rate analysis illustrated their efficiency, and the numerical results of convex optimization over fixed point sets demonstrated their effectiveness.
\end{minipage}
 \\[5mm]


\hbox to14cm{\hrulefill}\par

%%%%%%%%%%%%%%%%%%%%%%%%%%%%%%%%%%%%%%%%%%%%%%%%%%%%%%%%%%%%%%%%%%%%%%%
%%%%%%%%%%%%%
\section{Introduction}
\label{sec:1}
Stochastic programming problems have been recognized as significant, interesting problems that arise from practical applications in engineering and operational research. Stochastic optimization methods have thus been developed to efficiently solve various stochastic programming problems.

This paper considers a convex stochastic programming problem for which the objective function is given by the sum of convex functions or by a form of the expectation of convex functions and surveys stochastic optimization methods for solving it and related work. Incremental proximal point algorithms with randomized order \cite{bacak2014,bert2011} minimize the sum of convex functions. Random gradient and subgradient algorithms \cite{nedic2011} solve the problem of minimizing one convex function over sublevel sets of convex functions. It also discusses the connection between stochastic gradient descent and the randomized Kaczmarz method \cite{need2016}. Stochastic approximation and sample average approximation methods \cite{nemi2009} optimize the expected value of objective functions over a closed convex set. Incremental stochastic subgradient algorithms \cite{sun2009} minimize the sum of convex functions over a closed convex set. Stochastic approximation algorithms \cite{gha2012,gha2013,lan2012} perform convex stochastic composite optimization over a closed convex set. A distributed random projection algorithm \cite{lee2013} minimizes the sum of smooth, convex functions over the intersection of closed convex sets while incremental constraint projection-proximal methods \cite{bert2015} can be used to minimize the expected value of nonsmooth, convex functions over the intersection of closed convex sets. Multi-stage stochastic programming has been discussed \cite{berk2005,dyer2014,shapiro2008,zhao2005}, and the results \cite{bastin2006,gha2013_1,gha2016,gha_to,royset2012,xu2015} can even be applied to nonconvex stochastic optimization over the whole space or certain convex constraints. A proposed stochastic forward-backward splitting algorithm \cite{com2016} can be used to find the zeros of monotone operators.

In contrast to the stochastic programming considered in previous reports, this paper discusses stochastic programming problems in which each of the constraints is the {\em fixed point set} of a certain nonexpansive mapping. Convex optimization with fixed point constraints in a real Hilbert space is interesting and important because it enables consideration of optimization problems with complicated constraint sets onto which metric projections cannot be easily calculated and because it has many practical applications \cite{com,com1999,iiduka_siopt2013,iiduka2016,slav,yamada,yamada2011}. Although convex optimization with fixed point constraints has been analyzed in the deterministic case \cite{com,iiduka_siopt2013,iiduka_mp2014,iiduka2016,oms2016,iiduka_hishinuma_siopt2014,yamada} and stochastic fixed point algorithms have been presented \cite[Subchapter 10.3]{borkar2008}, \cite{com2015}, there have been no reports on stochastic optimization methods for convex optimization with fixed point constraints.
 
This paper is the first to consider convex stochastic programming problems with fixed point constraints and to present stochastic optimization algorithms for solving them. After the mathematical preliminaries and main problem statement are presented, the smooth convex stochastic programming problem is discussed (Section \ref{sec:3}), and a stochastic optimization algorithm is proposed to solve it. This algorithm (Algorithm \ref{algo:1}) blends a stochastic gradient method \cite[Subchapter 10.2]{borkar2008}, \cite{lee2013,nedic2011,nedic2001,need2016,sun2009} with the Halpern fixed point algorithm \cite{halpern,wit}, which is a useful fixed point algorithm. Next, the nonsmooth convex stochastic programming problem is discussed (Section \ref{sec:4}), and an algorithm (Algorithm \ref{algo:2}) is presented that is based on the stochastic proximal point algorithm \cite{bacak2014,bert2011,bert2015} and the Halpern fixed point algorithm.

One contribution of this paper is to enable consideration of (nonsmooth) convex stochastic optimization over fixed point sets of nonexpansive mappings, in contrast to recent papers \cite{iiduka_siopt2013,iiduka_mp2014,iiduka2016,iiduka_ejor2016,iiduka_hishinuma_siopt2014} that discussed deterministic convex optimization over the fixed point sets of nonexpansive mappings. The previous algorithm \cite{iiduka_mp2014} is a centralized acceleration algorithm for minimizing one smooth, strongly convex function over the fixed point set of a nonexpansive mapping. Although the algorithms in \cite{iiduka_siopt2013,iiduka_hishinuma_siopt2014} are decentralized algorithms that optimize the sum of smooth, convex objective functions over fixed point sets of nonexpansive mappings, they can work under the restricted situation such that the gradients of the functions are Lipschitz continuous and strongly or strictly monotone. The decentralized algorithms in \cite{iiduka2016,iiduka_ejor2016} can be applied to nonsmooth convex optimization with fixed point constraints. However, since the algorithms in \cite{iiduka2016,iiduka_ejor2016}, as with the previous algorithms \cite{iiduka_siopt2013,iiduka_mp2014,iiduka_hishinuma_siopt2014}, can be applied only to deterministic optimization, they cannot work on convex stochastic optimization over fixed point sets of nonexpansive mappings.

Another contribution is convergence analysis of the two proposed algorithms with diminishing step-size sequences. From the fact that a mapping containing the gradient of a convex function satisfies the nonexpansivity condition (Proposition \ref{nonexp}), it is shown that, under certain assumptions, any weak sequential cluster point of the sequence generated by the proposed gradient algorithm almost surely belongs to the solution set of the problem (Theorem \ref{thm:1}). The nonexpansivity condition of the proximity operator (Proposition \ref{prop:1}) means that, under certain assumptions, any weak sequential cluster point of the sequence generated by the proposed proximal point algorithm almost surely belongs to the solution set of the problem (Theorem \ref{thm:3}). Convergence rate analysis of the two algorithms is also provided (Propositions \ref{thm:2} and \ref{thm:4}).

This paper is organized as follows. Section \ref{sec:2} gives the mathematical preliminaries and states the main problem. Section \ref{sec:3} presents convergence and convergence rate analyses of the proposed gradient algorithm under certain assumptions. Section \ref{sec:4} presents convergence and convergence rate analyses of the proposed proximal point algorithm under certain assumptions. Section \ref{sec:5} considers specific convex optimization problems and compares numerically the behaviors of the two algorithms. Section \ref{sec:6} concludes the paper with a brief summary.
 
\section{Preliminaries}
\label{sec:2}
\subsection{Notation and definitions}
\label{subsec:2.1}
Let $H$ be a separable real Hilbert space with inner product $\langle \cdot,\cdot \rangle$, its induced norm $\| \cdot \|$, and Borel $\sigma$-algebra $\mathcal{B}$. Let $\mathbb{N}$ be the set of all positive integers including zero. Let $\mathrm{Id}$ denote the identity mapping on $H$. Let $\mathrm{Fix}(T) := \{x \in H \colon T(x) = x\}$ be the {\em fixed point set} of a mapping $T \colon H \to H$. The set of {\em weak sequential cluster points} \cite[Subchapters 1.7 and 2.5]{b-c} of a sequence $(x_n)_{n\in \mathbb{N}}$ in $H$ is denoted by $\mathcal{W}(x_n)_{n\in \mathbb{N}}$; i.e., $x\in \mathcal{W}(x_n)_{n\in\mathbb{N}}$ if and only if there exists a subsequence $(x_{n_l})_{l\in \mathbb{N}}$ of $(x_n)_{n\in \mathbb{N}}$ such that $(x_{n_l})_{l\in \mathbb{N}}$ weakly converges to $x$. Let $\mathbb{E}[X]$ denote the expectation of a random variable $X$. Given a probability space $(\Omega, \mathcal{F}, \mathbb{P})$, a $H$-valued random variable $x$ is defined by a measurable mapping $x \colon (\Omega,\mathcal{F}) \to (H,\mathcal{B})$. The $\sigma$-algebra generated by a family $\Phi$ of random variables is denoted by $\sigma(\Phi)$. Suppose that $(x_n)_{n\in \mathbb{N}}$ is a sequence of $H$-valued random variables and $C \subset H$. Then, any weak sequential cluster point of $(x_n)_{n\in \mathbb{N}}$ is said to almost surely belong to $C$ if there exists $\Omega_0 \in \mathcal{F}$ such that $\mathbb{P}(\Omega_0) = 1$ and $\mathcal{W}(x_n(\omega))_{n\in \mathbb{N}} \subset C$ for all $\omega \in \Omega_0$ (see also the proof of \cite[Corollary 2.7(i)]{com2015}). Suppose that $(x_n)_{n\in \mathbb{N}}$ and $(y_n)_{n\in \mathbb{N}}$ are positive real sequences. Let $O$ and $o$ denote Landau's symbols; i.e., $y_n = O(x_n)$ if there exist $c > 0$ and $n_0 \in \mathbb{N}$ such that $y_n \leq c x_n$ for all $n\geq n_0$, and $y_n = o(x_n)$ if, for all $\epsilon > 0$, there exists $n_0 \in \mathbb{N}$ such that $y_n \leq \epsilon x_n$ for all $n\geq n_0$.

A mapping $T \colon H \to H$ is said to be {\em nonexpansive} \cite[Definition 4.1(ii)]{b-c} if it is Lipschitz continuous with constant $1$; i.e., $\|T(x)-T(y)\| \leq \|x-y\|$ for all $x,y\in H$. $T$ is {\em firmly nonexpansive} \cite[Definition 4.1(i)]{b-c} if $\|T(x)-T(y)\|^2 + \|(\mathrm{Id}-T)(x) - (\mathrm{Id}-T)(y) \|^2 \leq \|x-y\|^2$ for all $x,y\in H$. This firm nonexpansivity condition obviously implies nonexpansivity. Given a nonempty, closed convex set $C \subset H$, the metric projection onto $C$, denoted by $P_C$, is defined for all $x\in H$ by $P_C(x) \in C$ and $\|x - P_C(x)\| = \inf_{y\in C} \|x-y\|$.

The {\em subdifferential} \cite[Definition 16.1, Corollary 16.14]{b-c} of a continuous, convex function $f \colon H \to \mathbb{R}$ is the set-valued operator $\partial f$ defined for all $x\in H$ by $\partial f (x) = \{ u\in H \colon f(y) \geq f(x) + \langle y-x,u \rangle \text{ } (y\in H) \} \neq \emptyset$. The condition $\partial f(x) = \{ \nabla f (x) \}$ holds for all $x\in H$ when $f$ is G\^{a}teaux differentiable \cite[Proposition 17.26]{b-c}. The {\em proximity operator} of $f$ \cite[Definition 12.23]{b-c}, \cite{minty1965,moreau1962}, denoted by $\mathrm{Prox}_f$, maps every $x\in H$ to the unique minimizer of $f(\cdot)+ (1/2) \| x - \cdot \|^2$.

\subsection{Main problem and propositions}
\label{subsec:2.2}
The following problem is considered in this paper.

\begin{prob}\label{prob:1}
Assume that 
\begin{enumerate}
\item[{\em (A1)}] $T^{(i)} \colon H \to H$ $(i\in \mathcal{I} := \{1,2,\ldots,I\})$ is firmly nonexpansive; 
\item[{\em (A2)}] $f^{(i)} \colon H \to \mathbb{R}$ $(i\in \mathcal{I})$ is convex and continuous.
\end{enumerate}
Then, our objective is to 
\begin{align*}
\text{minimize } f (x) := \mathbb{E}\left[ f^{(w)}(x) \right]
\text{ subject to } x \in X := \bigcap_{i\in \mathcal{I}} \mathrm{Fix}\left(T^{(i)}\right),
\end{align*}
where $f^{(w)}$ is a function involving a random variable $w \in \mathcal{I}$, and one assumes that 
\begin{enumerate}
\item[{\em (i)}] the solution set of the problem is nonempty;
\item[{\em (ii)}] there is an independent identically distributed sample $w_0, w_1, \ldots$ of realizations of the random variable $w$;
\item[{\em (iii)}] there is an oracle such that
\begin{itemize}
\item
for $(x,w) \in H \times \mathcal{I}$, it returns a stochastic firmly nonexpansive mapping $\mathsf{T}^{(w)}(x) := T^{(w)}(x)$;
\item
for $(z,w) \in H \times \mathcal{I}$, it returns a stochastic subgradient $\mathsf{G}^{(w)}(z) \in \partial f^{(w)}(z)$ or a stochastic proximal point $\mathsf{Prox}_{f^{(w)}}(z)$.
\end{itemize}
\end{enumerate}
\end{prob}

Problem \ref{prob:1} is discussed for the situation in which $(\mathsf{T}^{(w_n)}, f^{(w_n)})$ $(w_n \in \mathcal{I})$ is sampled at each iteration $n$. Let $J$ be the number of $f^{(i)}$. Even if $I < J$ (resp. $I > J$), the setting that $T^{(i)} := \mathrm{Id}$ $(i = I+1, I+2, \ldots,J)$ (resp. $f^{(j)}(x) := 0$ $(x\in H, j = J+1, J+2, \ldots,I)$), which satisfies (A1) (resp. (A2)), enables one to regard the stochastic optimization problem even when $J \neq I$ as Problem \ref{prob:1}.

The following propositions are used to prove the main theorems.

\begin{prop}{\em \cite[Proposition 2.3]{iiduka_JOTA}}\label{nonexp}
Let $f \colon H \to \mathbb{R}$ be convex and Fr\'echet differentiable, and let $\nabla f \colon H \to H$ be Lipschitz continuous with Lipschitz constant $L$. Then, $\mathrm{Id} - \lambda \nabla f$ is nonexpansive for all $\lambda \in [0,2/L]$.
\end{prop}

\begin{prop}{\em \cite[Propositions 12.26, 12.27, and 16.14]{b-c}}\label{prop:1}
Let $f\colon H \to \mathbb{R}$ be convex and continuous. Then, the following hold:
\begin{enumerate}
\item[{\em(i)}] Let $x,p\in H$. Then, $p = \mathrm{Prox}_f (x)$ if and only if $x-p \in \partial f (p)$. 
\item[{\em(ii)}] $\mathrm{Prox}_f$ is firmly nonexpansive with $\mathrm{Fix}(\mathrm{Prox}_f) = \argmin_{x\in H}f(x)$.
\item[{\em (iii)}] There exists $\delta > 0$ such that $\partial f(B(x;\delta))$ is bounded, where $B(x;\delta)$ represents a closed ball with center $x$ and radius $\delta$. 
\end{enumerate}
\end{prop}

\section{Stochastic gradient algorithm for smooth convex optimization}
\label{sec:3}
This section provides convergence properties of the following algorithm for solving Problem \ref{prob:1} when $f^{(i)}$ $(i\in \mathcal{I})$ is Fr\'echet differentiable. 

\begin{algorithm} 
\caption{Stochastic gradient algorithm for Problem \ref{prob:1}} 
\label{algo:1} 
\begin{algorithmic}[1] 
\REQUIRE
$n\in \mathbb{N}$, $(\alpha_n)_{n\in\mathbb{N}}, (\lambda_n)_{n\in \mathbb{N}} \subset (0,\infty)$.
\STATE
$n \gets 0$, 
$x_0 \in H$
\LOOP 
 \STATE 
 $y_{n} := \mathsf{T}^{(w_n)} \left(x_n - \lambda_n \mathsf{G}^{(w_n)}(x_n) \right)$
 \STATE 
 $x_{n+1} := \alpha_n x_0 + (1-\alpha_n) y_n$ 
 \STATE
 $n \gets n+1$
\ENDLOOP
\end{algorithmic}
\end{algorithm}

Algorithm \ref{algo:1} is obtained by blending the stochastic gradient method \cite[Subchapter 10.2]{borkar2008}, \cite{lee2013,nedic2011,need2016,sun2009} (i.e., $x_{n+1} = x_n -\lambda_n \mathsf{G}^{(w_n)}(x_n)$) with the Halpern fixed point algorithm \cite{halpern,wit}. The Halpern fixed point algorithm is defined by $x_0 \in H$ and $x_{n+1} = \alpha_n x_0 + (1-\alpha_n) T^{(i)}(x_n)$ $(n\in\mathbb{N})$ and converges strongly to a fixed point of $T^{(i)}$ when $(\alpha_n)_{n\in\mathbb{N}} \subset (0,1)$ satisfies $\lim_{n\to\infty} \alpha_n = 0$ and $\sum_{n=0}^\infty \alpha_n = \infty$. For Algorithm \ref{algo:1} to not only converge to a fixed point of $T^{(i)}$ but also to optimize $f^{(i)}$, Algorithm \ref{algo:1} needs to use an $(\alpha_n)_{n\in\mathbb{N}}$ that satisfies stronger conditions than $\lim_{n\to\infty} \alpha_n = 0$ and $\sum_{n=0}^\infty \alpha_n = \infty$ (see Assumption \ref{stepsize} for the conditions of $(\alpha_n)_{n\in\mathbb{N}}$ and $(\lambda_n)_{n\in\mathbb{N}}$).
 
\subsection{Assumptions for convergence analysis of Algorithm \ref{algo:1}}
\label{subsec:3.1}
Let us consider Problem \ref{prob:1} under (A1), (A2), and (A3) defined as follows. 
\begin{enumerate}
\item[(A3)]
$f^{(i)} \colon H \to \mathbb{R}$ $(i\in\mathcal{I})$ is Fr\'echet differentiable, and $\nabla f^{(i)} \colon H \to H$ is Lipschitz continuous with constant $L^{(i)}$.
\end{enumerate}

The following assumption is made.

\begin{assum}\label{stepsize}
Let $\sigma \geq 1$. The step-size sequences $(\alpha_n)_{n\in \mathbb{N}} \subset (0,1)$ and $(\lambda_n)_{n\in \mathbb{N}} \subset (0,1)$, which are monotone decreasing and converge to $0$, satisfy the following conditions:
\begin{align*}
&\text{{\em (C1)}} \sum_{n=0}^\infty \alpha_n = \infty, \text{ }
\text{{\em (C2)}} \lim_{n\to\infty} \frac{1}{\alpha_{n+1}} \left| \frac{1}{\lambda_{n+1}} - \frac{1}{\lambda_n} \right| = 0, \text{ } 
\text{{\em (C3)}} \lim_{n\to\infty} \frac{1}{\lambda_{n+1}} \left| 1 - \frac{\alpha_n}{\alpha_{n+1}} \right| = 0,\\
&\text{{\em (C4)}} \lim_{n\to\infty} \frac{\alpha_n}{\lambda_n} = 0, \text{ }
\text{{\em (C5)}} \frac{\alpha_n}{\alpha_{n+1}}, \frac{\lambda_n}{\lambda_{n+1}} \leq \sigma
\text{ } (n\in \mathbb{N}). 
\end{align*}
\end{assum}
Examples of $(\alpha_n)_{n\in \mathbb{N}}$ and $(\lambda_n)_{n\in \mathbb{N}}$ satisfying Assumption \ref{stepsize} are $\lambda_n := 1/(n+1)^a$ and $\alpha_n := 1/(n+1)^b$ $(n\in \mathbb{N})$, where $a\in (0,1/2)$ and $b\in (a,1-a)$. 

The collection of random variables is defined for all $n\in \mathbb{N}\backslash\{0\}$ by
\begin{align}\label{collection}
\mathcal{F}_n := 
\sigma(w_0, w_1, \ldots, w_{n-1}, y_0, y_1, \ldots, y_{n-1},x_0,x_1,\ldots,x_n).
\end{align}
Hence, given $\mathcal{F}_n$ defined by \eqref{collection}, the collection $y_0, y_1, \ldots, y_{n-1}$ and $x_0, x_1, \ldots, x_n$ generated by Algorithm \ref{algo:1} is determined.

The following is assumed for analyzing Algorithm \ref{algo:1}. 

\begin{assum}\label{omega}
The sequence $(w_n)_{n\in \mathbb{N}}$ satisfies the following conditions:
\begin{enumerate}
\item[{\em (i)}] For all $n\in \mathbb{N}$, there exists $m(n) \in \mathbb{N}$ such that $\bar{m} := \limsup_{n\to\infty} m(n) < \infty$ and $w_n = w_{n+m(n)}$ almost surely.
\item[{\em (ii)}] {\em \cite[Section 4 (see also Assumptions 4--7)]{bert2015}} There exists $\beta > 0$ such that, for all $i\in \mathcal{I}$ and for all $n\in \mathbb{N}$, $\beta \|x_n - T^{(i)}(x_n)\|^2 \leq \mathbb{E}[\|x_n - \mathsf{T}^{(w_n)}(x_n)\|^2 | \mathcal{F}_n]$ almost surely.
\end{enumerate}
Moreover, one of the following conditions holds.
\begin{enumerate}
\item[{\em (iii)}] {\em \cite[Section 5, Assumption 8]{bert2015}} $\mathbb{E}[f^{(w_n)}(x)|\mathcal{F}_n] = f(x)$ for all $x\in H$ and for all $n\in \mathbb{N}$ almost surely.
\item[{\em (iv)}] {\em \cite[Section 5, Assumption 9]{bert2015}} $(1/m) \sum_{l=tm}^{(t+1)m-1} \mathbb{E}[f^{(w_{l})}(x) | \mathcal{F}_{tm}] = f(x)$ for all $x\in H$ and for all $t\in \mathbb{N}$ almost surely, and $(\alpha_n)_{n\in \mathbb{N}}$ and $(\lambda_n)_{n\in \mathbb{N}}$ are constant within each cycle; i.e., $\alpha_{tm} = \alpha_{tm+1} = \cdots = \alpha_{(t+1)m-1}$ and $\lambda_{tm} = \lambda_{tm+1} = \cdots = \lambda_{(t+1)m-1}$.
\end{enumerate}
\end{assum}

A particularly interesting example of Sub-assumptions \ref{omega}(i) and (ii) is that, for all $t \in \mathbb{N}$, $(\mathsf{T}^{(w_n)})_{n\in \mathbb{N}}$, where $n = tI, tI +1, \ldots, (t+1)I-1$, is a permutation of $\{T^{(1)}, T^{(2)}, \ldots, T^{(I)}\}$ (see \cite[Subsection 4.3, Assumption 6]{bert2015} for the case in which $\mathsf{T}^{(w_n)}$ is a metric projection onto a simple, closed convex set).\footnote{Since all $T^{(i)}$ will be visited at least once within a cycle of $I$ iterations, Sub-assumption \ref{omega}(i) holds. From the nonexpansivity condition of a metric projection, the conclusions in \cite[Subsection 4.3]{bert2015} show that the sequence $(w_n)_{n\in\mathbb{N}}$ satisfies Sub-assumption \ref{omega}(ii) (see also Section \ref{sec:5}).} This enables one to consider the case in which the nonexpansive mappings are sampled in a cyclic manner (random shuffling or deterministic cycling). See Conditions (I), (II), and (IV) in Section \ref{sec:5} for other examples of $(w_n)_{n\in\mathbb{N}}$ satisfying Sub-assumptions \ref{omega}(i) and (ii).

Consider Problem \ref{prob:1} when $T := T^{(i)}$ $(i\in \mathcal{I})$ satisfying Sub-assumption \ref{omega}(ii), i.e.,
\begin{align}\label{prob:1_1}
\text{minimize } f(x) := \mathbb{E}\left[ f^{(w)}(x) \right] \text{ subject to } x\in \mathrm{Fix}(T).
\end{align}
Problem \eqref{prob:1_1} includes convex stochastic optimization problems in classifier ensemble \cite{hayashi,yin1,yin2}. However, the existing approaches in \cite{hayashi,yin1,yin2} are based on deterministic convex optimization and have not yet led to a complete solution of the classifier ensemble problem. Meanwhile, Theorem \ref{thm:1} guarantees that Algorithm \ref{algo:1} with $T := T^{(i)}$ $(i\in \mathcal{I})$, 
\begin{align}\label{algo:1_1}
x_{n+1} := \alpha_n x_0 + (1-\alpha_n) T \left( x_n - \lambda_n \mathsf{G}^{(w_n)}(x_n)\right) \text{ } (n\in \mathbb{N}),
\end{align}
can solve Problem \eqref{prob:1_1} including the classifier ensemble problem (see Subsection \ref{subsec:3.2} for convergence analysis of Algorithm \eqref{algo:1_1}).

Sub-assumption \ref{omega}(iii) implies that the sample component functions are conditionally unbiased \cite[Subsection 5.1, Assumption 8]{bert2015} while Sub-assumption \ref{omega}(iv) means that the functions are cyclically sampled \cite[Subsection 5.2, Assumption 9]{bert2015}. For simplicity, let us consider the case in which $(T^{(i)},f^{(i)})$ is sampled in a deterministic cyclic order (e.g., $w_0 = w_{tI} = I$, $w_{tI+i} = i$ $(t\in \mathbb{N}, i\in \mathcal{I})$). Then, Sub-assumption \ref{omega}(iv) means that $f(x) = (1/I) \sum_{i\in \mathcal{I}} f^{(i)}(x)$ $(x\in H)$. Problem \ref{prob:1} in such a deterministic case has been previously considered \cite{iiduka_siopt2013,iiduka_mp2014,iiduka2016,iiduka_ejor2016,iiduka_hishinuma_siopt2014}. In contrast to this deterministic case, Sub-assumptions \ref{omega}(i), (ii), and (iv) enable one to consider, for example, the stochastic Problem \ref{prob:1} with $f(x) = (1/I)\sum_{l=tI}^{(t+1)I-1}\mathbb{E}[f^{(w_{l})}(x)|\mathcal{F}_{tI}]$ $(x\in H, t\in \mathbb{N})$ for the case in which, for all $t \in \mathbb{N}$ and for a fixed $i_0 \in \mathcal{I}$, $(\mathsf{T}^{(w_n)}, f^{(w_n)})$ ($n=tI, tI+1, \ldots, (t+1)I-1$, $w_0 = w_{kI} = i_0$ $(k\in \mathbb{N})$) is sampled in a random cyclic order that differs depending on $t$. Section \ref{sec:5} provides numerical comparisons for the behaviors of Algorithm \ref{algo:1} with $(w_n)_{n\in\mathbb{N}}$ satisfying Assumption \ref{omega} (see (I)--(IV) in Section \ref{sec:5}).

The convergence of Algorithm \ref{algo:1} depends on the following assumption.
\begin{assum}\label{bounded}
The sequence $(y_n)_{n\in \mathbb{N}}$ is almost surely bounded. 
\end{assum}

Assumption \ref{bounded} and the definition of $(x_n)_{n\in \mathbb{N}}$ ensure that $(x_n)_{n\in \mathbb{N}}$ is almost surely bounded. This guarantees that there exist $\bar{\Omega} \in \mathcal{F}$ with $\mathbb{P}(\bar{\Omega}) = 1$ and a weak sequential cluster point of $(x_n(\omega))_{n\in\mathbb{N}}$ $(\omega \in \bar{\Omega})$ in Algorithm \ref{algo:1}; i.e., there exists a weak convergent subsequence $(x_{n_i}(\omega))_{i\in \mathbb{N}}$ of $(x_n(\omega))_{n\in \mathbb{N}}$ $(\omega \in \bar{\Omega})$. Hence, Assumption \ref{bounded} is needed to analyze the weak convergence of Algorithm \ref{algo:1}. Suppose that a bounded, closed convex set $C \subset H$ can be chosen in advance such that the metric projection onto $C\supset X$, denoted by $P_C$, is easily computed within a finite number of arithmetic operations \cite[Subchapter 28]{b-c} (e.g., $C$ is a closed ball with a large enough radius). Then, $y_n$ $(n\in \mathbb{N})$ in Algorithm \ref{algo:1} (step 3 in Algorithm \ref{algo:1}) can be replaced with 
\begin{align}\label{y_n}
y_n := P_C \left[\mathsf{T}^{(w_n)} \left(x_n - \lambda_n \mathsf{G}^{(w_n)}(x_n) \right) \right],
\end{align}
which means that Assumption \ref{bounded} holds. The same discussion in Subsection \ref{subsec:3.2} ensures that any weak sequential cluster point of the sequence $(x_n)_{n\in \mathbb{N}}$ generated by Algorithm \ref{algo:1} with \eqref{y_n} belongs to the solution set of Problem \ref{prob:1} without assuming Assumption \ref{bounded}.

\subsection{Convergence analysis of Algorithm \ref{algo:1}}
\label{subsec:3.2}
The convergence of Algorithm \ref{algo:1} can be analyzed as follows.

\begin{thm}\label{thm:1}
Suppose that Assumptions (A1)-(A3) and \ref{stepsize}-\ref{bounded} hold, and let $(x_n)_{n\in \mathbb{N}}$ be the sequence generated by Algorithm \ref{algo:1}. Then, any weak sequential cluster point of $(x_n)_{n\in\mathbb{N}}$ almost surely belongs to the solution set of Problem \ref{prob:1}.
\end{thm}

The proof of Theorem \ref{thm:1} is divided into five steps (Lemmas \ref{lem:1}, \ref{lem:2}, \ref{lem:3}, \ref{lem:4}, and the proof of Theorem \ref{thm:1}). First, the following lemma is proven.

\begin{lem}\label{lem:1}
Suppose that Assumptions (A1)-(A3), \ref{stepsize}, \ref{omega}(i), and \ref{bounded} hold. Then, almost surely
\begin{align*}
\lim_{n\to\infty} \mathbb{E} \left[\frac{\| x_{n+m+1} - x_{n+1} \|}{\lambda_{n+m}} 
\bigg| \mathcal{F}_n \right] = 0.
\end{align*} 
\end{lem}

{\em Proof:}
Assumption \ref{bounded} means the almost sure boundedness of $(x_n)_{n\in \mathbb{N}}$. Accordingly, the Lipschitz continuity of $\nabla f^{(i)}$ $(i\in \mathcal{I})$ (see (A3)) leads to the almost sure boundedness of $(\nabla f^{(i)}(x_n))_{n\in \mathbb{N}}$ $(i\in \mathcal{I})$; i.e., $M_1 := \max_{i\in \mathcal{I}} \{\sup_{n\in \mathbb{N}} \| \nabla f^{(i)} (x_n)\| \}< \infty$ almost surely. From the monotone decreasing condition of $(\lambda_n)_{n\in \mathbb{N}}$, there exists $n_0\in \mathbb{N}$ such that, for all $n \geq n_0$, $\lambda_n \leq L:= 2/\max_{i\in \mathcal{I}} L^{(i)}$. Hence, (A2), (A3), and Proposition \ref{nonexp} imply that $\mathrm{Id} - \lambda_n \nabla f^{(i)}$ $(i\in \mathcal{I}, n\geq n_0)$ is nonexpansive. Sub-assumption \ref{omega}(i) ensures that, for all $n\geq n_0$, there exists $m(n) \in \mathbb{N}$ such that $\limsup_{n\to \infty} m(n) < \infty$, $\mathsf{T}^{(w_{n+m})} = \mathsf{T}^{(w_n)}$, and $f^{(w_{n+m})} = f^{(w_n)}$ almost surely. Accordingly, (A1) and the triangle inequality ensure that, for all $n\geq n_0$, almost surely
\begin{align*}
\left\| y_{n+m} - y_n \right\|
= &\left\| \mathsf{T}^{(w_n)} \left(x_{n+m} - \lambda_{n+m} \mathsf{G}^{(w_{n})}(x_{n+m}) \right) 
 - \mathsf{T}^{(w_n)} \left(x_n - \lambda_n \mathsf{G}^{(w_n)}(x_n) \right)\right\|\\
\leq & \left\| \left(x_{n+m} - \lambda_{n+m} \mathsf{G}^{(w_{n})}(x_{n+m}) \right) 
 - \left(x_n - \lambda_n \mathsf{G}^{(w_n)}(x_n) \right)\right\|\\
\leq & \left\| \left(x_{n+m} - \lambda_{n+m} \mathsf{G}^{(w_{n})}(x_{n+m}) \right) 
 - \left(x_n - \lambda_{n+m} \mathsf{G}^{(w_{n})}(x_n) \right)\right\|\\
 &+ \left| \lambda_{n+m} - \lambda_n \right| \left\| \mathsf{G}^{(w_n)}(x_n) \right\|, 
\end{align*}
which, together with the nonexpansivity of $\mathrm{Id} - \lambda_{n+m} \mathsf{G}^{(w_{n})}$, implies that
\begin{align*}
\left\| y_{n+m} - y_n \right\|
\leq \left\| x_{n+m} - x_n \right\| + M_1 \left| \lambda_{n+m} - \lambda_n \right|.
\end{align*}
Since the definition of $x_n$ $(n\in \mathbb{N})$ and the triangle inequality mean that, for all $n\geq n_0$, 
\begin{align*}
\| x_{n+m+1} - x_{n+1} \|
&= \left\| \left( \alpha_{n+m} - \alpha_n \right) \left(x_0 - y_n\right) + \left(1 - \alpha_{n+m} \right) \left(y_{n+m} - y_n \right) \right\|\\
&\leq \left(1 - \alpha_{n+m} \right) \left\| y_{n+m} - y_n \right\| 
 + \left| \alpha_{n+m} - \alpha_n \right| \left\|x_0 - y_n \right\|, 
\end{align*}
meaning that, for all $n\geq n_0$, almost surely
\begin{align}\label{xn}
\begin{split}
\left\| x_{n+m+1} - x_{n+1} \right\| 
&\leq \left(1 - \alpha_{n+m} \right) \left\{ \left\| x_{n+m} - x_n \right\| + M_1 \left| \lambda_{n+m} - \lambda_n \right|\right\}\\ 
&\quad + \left| \alpha_{n+m} - \alpha_n \right| \left\|x_0 - y_n\right\|\\
&\leq \left(1 - \alpha_{n+m} \right)\left\| x_{n+m} - x_n \right\| + M_1 \left| \lambda_{n+m} - \lambda_n \right|\\
&\quad + M_2 \left| \alpha_{n+m} - \alpha_n \right|,
\end{split}
\end{align}
where almost surely $M_2 := \sup_{n\in \mathbb{N}} \|y_n - x_0\| < \infty$. Therefore, for all $n\geq n_0$, almost surely
\begin{align*}
&\frac{\left\| x_{n+m+1} - x_{n+1} \right\|}{\lambda_{n+m}}\\ 
\leq& 
\left(1 - \alpha_{n+m} \right) \frac{\left\| x_{n+m} - x_n \right\|}{\lambda_{n+m}} 
 + M_1 \frac{\left| \lambda_{n+m} - \lambda_n \right|}{\lambda_{n+m}}
 + M_2 \frac{\left| \alpha_{n+m} - \alpha_n \right|}{\lambda_{n+m}}\\
=& 
\left(1 - \alpha_{n+m} \right) \frac{\left\| x_{n+m} - x_n \right\|}{\lambda_{n+m-1}}
 +\left(1 - \alpha_{n+m} \right) \left\{ \frac{\left\| x_{n+m} - x_n \right\|}{\lambda_{n+m}} - 
 \frac{\left\| x_{n+m} - x_n \right\|}{\lambda_{n+m-1}} \right\}\\
 &+ M_1 \frac{\left| \lambda_{n+m} - \lambda_n \right|}{\lambda_{n+m}}
 + M_2 \frac{\left| \alpha_{n+m} - \alpha_n \right|}{\lambda_{n+m}}\\
\leq&
\left(1 - \alpha_{n+m} \right) \frac{\left\| x_{n+m} - x_n \right\|}{\lambda_{n+m-1}}
 + M_3 \left| \frac{1}{\lambda_{n+m}} - \frac{1}{\lambda_{n+m-1}} \right|
 + M_1 \frac{\left| \lambda_{n+m} - \lambda_n \right|}{\lambda_{n+m}}\\
 &+ M_2 \frac{\left| \alpha_{n+m} - \alpha_n \right|}{\lambda_{n+m}},
\end{align*}
where almost surely $M_3 := \sup_{n\in \mathbb{N}} \| x_{n+m} - x_n \| < \infty$. Accordingly, for all $n\geq n_0$, almost surely
\begin{align}
\frac{\left\| x_{n+m+1} - x_{n+1} \right\|}{\lambda_{n+m}}
\leq&
\left(1 - \alpha_{n+m} \right) \frac{\left\| x_{n+m} - x_n \right\|}{\lambda_{n+m-1}}
 + L \alpha_{n+m} \frac{M_1}{\alpha_{n+m}} \left| \frac{1}{\lambda_n} - \frac{1}{\lambda_{n+m}} \right|\label{ineq:1}\\
 &+ \alpha_{n+m} \frac{M_3}{\alpha_{n+m}} \left| \frac{1}{\lambda_{n+m}} - \frac{1}{\lambda_{n+m-1}} \right|
 + \alpha_{n+m} \frac{M_2}{\lambda_{n+m}} \left| 1 - \frac{\alpha_n}{\alpha_{n+m}} \right|,\nonumber
\end{align}
where the second term on the right comes from $\lambda_n \leq L$ $(n\geq n_0)$ and 
\begin{align*}
\frac{\left| \lambda_{n+m} - \lambda_n \right|}{\lambda_{n+m}}
=
L \frac{\left| \lambda_{n+m} - \lambda_n \right|}{L \lambda_{n+m}}
\leq L \frac{\left| \lambda_{n+m} - \lambda_n \right|}{\lambda_{n} \lambda_{n+m}}
= L \left| \frac{1}{\lambda_n} - \frac{1}{\lambda_{n+m}} \right|.
\end{align*}
Condition (C5) and the triangle inequality mean that, for all $n \geq n_0$ and for all $l \geq 1$,
\begin{align*}
\frac{1}{\alpha_{n+l+1}} \left| \frac{1}{\lambda_n} - \frac{1}{\lambda_{n+l+1}} \right|
&\leq 
\frac{\alpha_{n+l}}{\alpha_{n+l+1}} \frac{1}{\alpha_{n+l}}
\left| \frac{1}{\lambda_n} - \frac{1}{\lambda_{n+l}} \right|
+ \frac{1}{\alpha_{n+l+1}}
\left| \frac{1}{\lambda_{n+l}} - \frac{1}{\lambda_{n+l+1}} \right|\\
&\leq 
\sigma \frac{1}{\alpha_{n+l}}
\left| \frac{1}{\lambda_n} - \frac{1}{\lambda_{n+l}} \right|
+ \frac{1}{\alpha_{n+l+1}}
\left| \frac{1}{\lambda_{n+l}} - \frac{1}{\lambda_{n+l+1}} \right|,\\
\frac{1}{\lambda_{n+l+1}} \left| 1 - \frac{\alpha_n}{\alpha_{n+l+1}} \right|
&\leq
\frac{\alpha_{n+l}}{\alpha_{n+l+1}} \frac{\lambda_{n+l}}{\lambda_{n+l+1}}
\frac{1}{\lambda_{n+l}} \left| 1 - \frac{\alpha_n}{\alpha_{n+l}} \right|
+
\frac{1}{\lambda_{n+l+1}} \left| 1 - \frac{\alpha_{n+l}}{\alpha_{n+l+1}} \right|\\
&\leq \sigma^2 \frac{1}{\lambda_{n+l}} \left| 1 - \frac{\alpha_n}{\alpha_{n+l}} \right|
+
\frac{1}{\lambda_{n+l+1}} \left| 1 - \frac{\alpha_{n+l}}{\alpha_{n+l+1}} \right|.
\end{align*}
Conditions (C2) and (C3) thus mean that, for all $l\geq 1$,
\begin{align}\label{step}
\lim_{n\to\infty} \frac{1}{\alpha_{n+l}}
\left| \frac{1}{\lambda_n} - \frac{1}{\lambda_{n+l}} \right| = 0
\text{ and } 
\lim_{n\to\infty} \frac{1}{\lambda_{n+l}} \left| 1 - \frac{\alpha_n}{\alpha_{n+l}} \right| = 0.
\end{align}
Hence, (C2) and \eqref{step} guarantee that, for all $\epsilon > 0$, there exists $n_1 \in \mathbb{N}$ such that, for all $n \geq n_1$,
\begin{align*}
\frac{M_1 L}{\alpha_{n+m}} \left| \frac{1}{\lambda_n} - \frac{1}{\lambda_{n+m}} \right| \leq \frac{\epsilon}{3}, \text{ } 
\frac{M_3}{\alpha_{n+m}} \left| \frac{1}{\lambda_{n+m}} - \frac{1}{\lambda_{n+m-1}} \right| \leq \frac{\epsilon}{3}, \text{ }
\frac{M_2}{\lambda_{n+m}} \left| 1 - \frac{\alpha_n}{\alpha_{n+m}} \right| \leq \frac{\epsilon}{3}.
\end{align*}
Therefore, \eqref{ineq:1} means that, for all $n \geq n_2 := \max \{ n_0, n_1 \}$, almost surely
\begin{align}\label{ine}
\frac{\left\| x_{n+m+1} - x_{n+1} \right\|}{\lambda_{n+m}}
\leq&
\left(1 - \alpha_{n+m} \right) \frac{\left\| x_{n+m} - x_n \right\|}{\lambda_{n+m-1}}
 + \epsilon \alpha_{n+m}.
\end{align}
Further, induction guarantees that, for all $n \geq n_2$, almost surely
\begin{align*}
&\quad \frac{\left\| x_{n+1+ m(n)} - x_{n+1} \right\|}{\lambda_{n+m(n)}}\\
&\leq
\left(1 - \alpha_{n+m(n)} \right) 
 \bigg\{ \left(1 - \alpha_{n-1+m(n-1)} \right) \frac{\left\| x_{n-1+m(n-1)} - x_{n-1} \right\|}{\lambda_{n-2+m(n-1)}}\\
&\quad + \epsilon \left(1- \left(1- \alpha_{n-1+m(n-1)}\right) \right)
 \bigg\}
 + \epsilon \alpha_{n+m(n)}\\
&= 
\left(1 - \alpha_{n+m(n)} \right) \left(1 - \alpha_{n-1+m(n-1)} \right)
 \frac{\left\| x_{n-1+m(n-1)} - x_{n-1} \right\|}{\lambda_{n-2+m(n-1)}}\\ 
 &\quad + \epsilon \left\{ 1 - \left(1 - \alpha_{n+m(n)} \right) \left(1 - \alpha_{n-1+m(n-1)} \right) \right\}\\
&\leq
\prod_{k=n_2}^{n} \left(1 - \alpha_{k+m(k)} \right) \frac{\left\| x_{n_2+m(n_2)} - x_{n_2} \right\|}{\lambda_{n_2 -1 +m(n_2)}}
 + \epsilon \left\{1 - \prod_{k=n_2}^{n} \left(1 - \alpha_{k+m(k)} \right) \right\}.
\end{align*}
By taking the expectation in this inequality conditioned on $\mathcal{F}_{n}$ $(n\geq n_2)$ defined in \eqref{collection}, we have for all $n\geq n_2$
\begin{align}
\mathbb{E} \left[\frac{\left\| x_{n+m(n)+1} - x_{n+1} \right\|}{\lambda_{n+m(n)}} \Bigg| \mathcal{F}_{n} \right]
&\leq
\prod_{k=n_2}^{n} \left(1 - \alpha_{k+m(k)} \right) \mathbb{E} \left[ \frac{\left\| x_{n_2+m(n_2)} - x_{n_2} \right\|}{\lambda_{n_2+m(n_2)-1}}\Bigg| \mathcal{F}_{n} \right]\nonumber\\
&\quad + \epsilon \left\{1 - \prod_{k=n_2}^{n} \left(1 - \alpha_{k+m(k)} \right) \right\}\label{ineq:2}
\end{align}
almost surely. Moreover, Sub-assumption \ref{omega}(i) means the existence of $\hat{m}\in \mathbb{N}$ satisfying $\max\{ m(k) \colon k=n,n-1,\ldots,n_2 \} \leq \hat{m}$. Accordingly, Condition (C1) and the monotone decreasing condition of $(\alpha_n)_{n\in \mathbb{N}}$ lead to the finding that $0\leq \limsup_{n\to\infty} \prod_{k=n_2}^{n} (1 - \alpha_{k+m(k)}) \leq \limsup_{n\to\infty} \prod_{k=n_2}^{n} (1 - \alpha_{k+\hat{m}} ) =0$. Therefore, \eqref{ineq:2} means that, almost surely
\begin{align*}
\limsup_{n\to\infty} \mathbb{E} \left[\frac{\left\| x_{n+m(n)+1} - x_{n+1} \right\|}{\lambda_{n+m(n)}} \Bigg| \mathcal{F}_n \right]
%&=
%\limsup_{l\to\infty} \mathrm{E} \left[\frac{\left\| x_{l+n_2 +m(l+n_2)+1} - x_{l+n_2+1} \right\|}{\lambda_{l+n_2+m(l+n_2)}} \Bigg| \mathcal{F}_{l+n_2} \right]\\
&\leq 
 \epsilon,
\end{align*}
which, together with the arbitrary condition of $\epsilon$, means that Lemma \ref{lem:1} holds.

Lemma \ref{lem:1} leads to the following.

\begin{lem}\label{lem:2}
Suppose that the assumptions in Lemma \ref{lem:1} hold. Then, almost surely
\begin{align*}
\lim_{n\to\infty} \mathbb{E} \left[\left\|x_{n} - y_{n} \right\|^2
\Big| \mathcal{F}_n \right] = 0 \text{ and }
\lim_{n\to\infty} \mathbb{E}\left[ \left\| x_{n} - \mathsf{T}^{(w_n)} (x_{n}) \right\|^2 \bigg| \mathcal{F}_n \right] = 0.
\end{align*}
\end{lem}

{\em Proof:}
Fix $x\in X \subset \mathrm{Fix}(T^{(i)})$ $(i\in \mathcal{I})$ and $n\in \mathbb{N}$ arbitrarily. Assumption (A1) ensures that, for all $k\in \mathbb{N}$, $\| y_k - x \|^2 \leq \| (x_k - x ) - \lambda_k \mathsf{G}^{(w_k)}(x_k)\|^2 - \| (x_k - y_k ) - \lambda_k \mathsf{G}^{(w_k)}(x_k) \|^2$. Hence, from $\| x -y\|^2 = \|x\|^2 -2 \langle x,y\rangle +\|y\|^2$ $(x,y\in H)$, 
\begin{align*}
\left\| y_k - x \right\|^2
%&= \left\| T^{(w_k)} \left(x_k - \lambda_k \nabla f^{(w_k)}(x_k) \right) - T^{(w_k)} (x) \right\|^2\\
%&\leq \left\| \left(x_k - x \right) - \lambda_k \nabla f^{(w_k)}(x_k) \right\|^2 - \left\| \left(x_k - y_k \right) - \lambda_k \nabla f^{(w_k)}(x_k)\right\|^2\\
&\leq \left\|x_k - x \right\|^2 + 2 \lambda_k \left\langle x - y_k, 
\mathsf{G}^{(w_k)}(x_k) \right\rangle
- \left\|x_k - y_k \right\|^2.
\end{align*}
The definition of $x_k$ $(k\in \mathbb{N})$ and the convexity of $\| \cdot \|^2$ thus imply that, for all $k\in \mathbb{N}$,
\begin{align*}
\left\| x_{k+1} - x \right\|^2
%&\leq
%\alpha_k \left\| x_0 - x \right\|^2 + (1-\alpha_k) \left\| y_k - x\right\|^2\\
%&\leq 
%\alpha_k \left\| x_0 - x \right\|^2 + (1-\alpha_k)
%\left\{ \left\|x_k - x \right\|^2 + 2 \lambda_k \left\langle x - y_k, \nabla f^{(w_k)}(x_k) \right\rangle
%- \left\|x_k - y_k \right\|^2 \right\}\\
&\leq 
\alpha_k \left\| x_0 - x \right\|^2 + \left\|x_k - x \right\|^2
+ 2(1-\alpha_k)\lambda_k \left\langle x - y_k, \mathsf{G}^{(w_k)}(x_k) \right\rangle\\
&\quad - (1-\alpha_k)\left\|x_k - y_k \right\|^2.
\end{align*}
Since the above inequality holds for $k = n+m(n), n+m(n)-1, \ldots, n+1$, it can be deduced that
\begin{align*}
\left\| x_{n+m+1} - x \right\|^2
&\leq 
 \left\|x_{n+1} - x \right\|^2 + \left\| x_0 - x \right\|^2
 \sum_{k=n+1}^{n+m} \alpha_{k}
 - \sum_{k=n+1}^{n+m} (1-\alpha_{k}) \left\|x_{k} - y_{k} \right\|^2\\
&\quad + 2 \sum_{k=n+1}^{n+m} \lambda_{k} \left|\left\langle x - y_{k}, \mathsf{G}^{(w_k)}(x_{k}) \right\rangle \right|,
\end{align*}
which, together with $M_4 := \sup_{n\in\mathbb{N}} 2|\langle x - y_{n},  \mathsf{G}^{(w_n)}(x_{n}) \rangle | < \infty$ almost surely, and the triangle inequality, means that, almost surely
\begin{align}\label{xnyn3}
(1-\alpha_{n+1}) \left\|x_{n+1} - y_{n+1} \right\|^2
%&\leq
%\left\| x_0 - x \right\|^2 \sum_{k=n+1}^{n+m} \alpha_{k}
%+ \left\|x_{n+1} - x \right\|^2 - \left\| x_{n+m+1} - x \right\|^2
%+ M_4 \sum_{k=n+1}^{n+m} \lambda_{k} \nonumber\\
%&= 
%\left\| x_0 - x \right\|^2 \sum_{k=n+1}^{n+m} \alpha_{k} 
%+ 2 M_1 \sum_{k=n+1}^{n+m} \lambda_{k} \left\|x - y_{k} \right\|\nonumber\\
%&\quad + \left(\left\|x_{n+1} - x \right\| 
% + \left\| x_{n+m+1} - x \right\| \right)
% \left(\left\|x_{n+1} - x \right\| 
% - \left\| x_{n+m+1} - x \right\| \right)\nonumber\\
&\leq 
\left\| x_0 - x \right\|^2 \sum_{k=n+1}^{n+m} \alpha_{k} 
+ M_4 \sum_{k=n+1}^{n+m} \lambda_{k} \nonumber\\ 
&\quad + \lambda_{n+m} \left(\left\|x_{n+1} - x \right\| 
 + \left\| x_{n+m+1} - x \right\| \right)
 \frac{\left\|x_{n+1} - x_{n+m+1} \right\|}{\lambda_{n+m}}. 
\end{align}
Taking the expectation in this inequality conditioned on $\mathcal{F}_{n+1}$ defined in \eqref{collection} leads to the finding that, almost surely
\begin{align}\label{xnyn2}
\begin{split}
&\quad (1-\alpha_{n+1}) \mathbb{E} \left[\left\|x_{n+1} - y_{n+1} \right\|^2
\Big| \mathcal{F}_{n+1} \right]\\
&\leq 
\left\| x_0 - x \right\|^2 \sum_{k=n+1}^{n+m} \alpha_{k} 
+ M_4 \sum_{k=n+1}^{n+m} \lambda_{k}\\
&\quad 
+ \lambda_{n+m} \mathbb{E} \left[ \left(\left\|x_{n+1} - x \right\| 
 + \left\| x_{n+m+1} - x \right\| \right)
 \frac{\left\|x_{n+1} - x_{n+m+1} \right\|}{\lambda_{n+m}} \bigg| \mathcal{F}_{n+1} \right].
\end{split} 
\end{align}
Hence, from the definition of $\mathcal{F}_n$ $(n\in\mathbb{N})$, Assumption \ref{bounded}, Lemma \ref{lem:1}, and $\lim_{n\to\infty} \alpha_n = \lim_{n\to \infty} \lambda_n = 0$, we have  
\begin{align}\label{xnyn}
\lim_{n\to\infty} \mathbb{E} \left[\left\|x_{n} - y_{n} \right\|^2
\Big| \mathcal{F}_n \right] = 0
\end{align}
almost surely. Further, since (A1) means that, for all $n\in \mathbb{N}$,
\begin{align*}
\left\| y_{n} - \mathsf{T}^{(w_n)} (x_{n}) \right\|
&=
\left\| \mathsf{T}^{(w_n)} \left(x_{n} - \lambda_{n} \mathsf{G}^{(w_n)}(x_{n}) \right) - \mathsf{T}^{(w_n)} (x_{n}) \right\|
\leq 
\lambda_{n} \left\| \mathsf{G}^{(w_n)}(x_{n}) \right\|,
\end{align*}
we find that, for all $n\in \mathbb{N}$,
\begin{align}\label{xnt}
\begin{split}
\left\| x_{n} - \mathsf{T}^{(w_n)} (x_{n}) \right\|^2
&\leq 
2 \left\| x_{n} - y_{n} \right\|^2 + 2 \left\| y_{n} - \mathsf{T}^{(w_n)} (x_{n})\right\|^2\\
&\leq 
2 \left\| x_{n} - y_{n} \right\|^2 + 2 \lambda_{n}^2 \left\| \mathsf{G}^{(w_n)}(x_{n}) \right\|^2,
\end{split}
\end{align}
where the first inequality comes from $\| x+ y \|^2 \leq 2 \|x\|^2 + 2 \|y\|^2$ $(x,y\in H)$. Accordingly, \eqref{xnyn}, Assumption \ref{bounded}, and the convergence of $(\lambda_n)_{n\in \mathbb{N}}$ to $0$ guarantee that, almost surely $\lim_{n\to\infty} \mathbb{E} [\| x_{n} - \mathsf{T}^{(w_n)} (x_{n})\|^2| \mathcal{F}_n] = 0$. This completes the proof.

The following lemma demonstrates that any weak sequential cluster point of $(x_n)_{n\in \mathbb{N}}$ in Algorithm \ref{algo:1} is almost surely in $X$. 

\begin{lem}\label{lem:3}
Suppose that Sub-assumption \ref{omega}(ii) and the assumptions in Lemma \ref{lem:1} hold. Then, for all $i\in \mathcal{I}$, almost surely
\begin{align*}
\lim_{n\to\infty} \left\|x_n -T^{(i)}(x_n) \right\| = 0
\text{ and } 
\lim_{n\to\infty} \left\|x_n -T^{(i)}\left(x_n - \lambda_n \nabla f^{(i)}(x_n)\right) \right\| = 0.
\end{align*}
\end{lem}

{\em Proof:}
Sub-assumption \ref{omega}(ii) and Lemma \ref{lem:2} guarantee that, for all $j\in \mathcal{I}$, almost surely
\begin{align*}
\beta \limsup_{n\to \infty} \left\| x_n - T^{(j)} (x_n) \right\|^2 
\leq \lim_{n\to\infty} \mathbb{E} \left[ \left\| x_n - \mathsf{T}^{(w_n)} (x_n) \right\|^2 \bigg| \mathcal{F}_n \right] = 0.
\end{align*}
This means that $\lim_{n\to\infty} \|x_n -T^{(j)}(x_n) \|$ $(j\in \mathcal{I})$ almost surely equals $0$. The triangle inequality and (A1) ensure that, for all $i\in \mathcal{I}$ and for all $n\in \mathbb{N}$,
\begin{align*}
\left\|x_n -T^{(i)}\left(x_n - \lambda_n \nabla f^{(i)}(x_n)\right) \right\|
\leq \left\|x_n -T^{(i)}(x_n)\right\|
+ \lambda_n \left\| \nabla f^{(i)}(x_n) \right\|,
\end{align*}
which, together with Assumption \ref{bounded}, $\lim_{n\to\infty}\lambda_n = 0$ almost surely, and $\lim_{n\to\infty}\|x_n - T^{(i)}(x_n)\|=0$ almost surely, means that $\lim_{n\to\infty} \|x_n -T^{(i)}(x_n - \lambda_n \nabla f^{(i)}(x_n))\|$ $(j\in \mathcal{I})$ almost surely equals $0$. This completes the proof.

The following can also be proved.

\begin{lem}\label{lem:4}
Suppose that the assumptions in Theorem \ref{thm:1} hold. Then, almost surely
\begin{align*}
\limsup_{n\to\infty} f(x_n) \leq f^\star := \min_{x\in X} f(x).
\end{align*}
\end{lem}

{\em Proof:}
Fix $x^\star \in X^\star := \{ x^\star \in X \colon f(x^\star) = f^\star \}$ and $n\in \mathbb{N}$ arbitrarily. From (A1), for all $k\in \mathbb{N}$, $\| y_k - x^\star \|^2 \leq \| (x_k - x^\star ) - \lambda_k \mathsf{G}^{(w_k)}(x_k) \|^2$, which, together with $\| x-y\|^2 = \|x\|^2 -2\langle x,y\rangle +\|y\|^2$ $(x,y\in H)$ and the definition of $\partial f$, means that, for all $k\in \mathbb{N}$, almost surely
\begin{align*}
\left\| y_k - x^\star \right\|^2
\leq \left\| x_k - x^\star \right\|^2 + 2 \lambda_k \left( f^{(w_k)}(x^\star) - f^{(w_k)}(x_k) \right) + M_1^2 \lambda_k^2.
\end{align*}
Hence, the convexity of $\|\cdot \|^2$ means that, for all $k\in \mathbb{N}$, almost surely
\begin{align*}
&\quad \left\| x_{k+1} - x^\star \right\|^2\\
%&\leq
%\alpha_k \left\| x_0 - x^\star \right\|^2 + (1-\alpha_k) \left\| y_k - x^\star\right\|^2\\
%&\leq 
%\alpha_k \left\| x_0 - x^\star \right\|^2 + (1-\alpha_k)
%\left\{ \left\| x_k - x^\star \right\|^2 + 2 \lambda_k \left( f^{(w_k)}(x^\star) - f^{(w_k)}(x_k) \right) + M_1^2 \lambda_k^2 \right\}\\
&\leq 
\alpha_k \left\| x_0 - x^\star \right\|^2 + \left\| x_k - x^\star \right\|^2
+ 2 (1-\alpha_k)\lambda_k \left( f^{(w_k)}(x^\star) - f^{(w_k)}(x_k) \right) + M_1^2 \lambda_k^2.
\end{align*}
Since the above inequality holds for $k = n+m(n), n+m(n)-1, \ldots, n+1$, almost surely
\begin{align}\label{ineq:f}
\begin{split}
&\quad \frac{2}{\lambda_{n+m}} \sum_{k=n+1}^{n+m} (1-\alpha_k)\lambda_k \left( f^{(w_k)}(x_k) - f^{(w_k)}(x^\star) \right)\\
&\leq M_5 \frac{\left\| x_{n+m+1} - x_{n+1} \right\|}{\lambda_{n+m}}
+ \frac{\left\| x_0 - x^\star \right\|^2}{\lambda_{n+m}} \sum_{k=n+1}^{n+m} \alpha_k
 + \frac{M_1^2}{\lambda_{n+m}} \sum_{k=n+1}^{n+m} \lambda_k^2,
\end{split} 
\end{align}
where almost surely $M_5 := \sup_{n\in \mathbb{N}}( \| x_{n+1} - x^\star \| + \| x_{n+m+1} - x^\star \| ) < \infty$. 

Now, let us assume that Sub-assumption \ref{omega}(iii) holds. Then, for all $x\in H$, almost surely $\mathbb{E}[f^{(w_{n+1})}(x)|\mathcal{F}_n] = \mathbb{E}[\mathbb{E}[f^{(w_{n+1})}(x)|\mathcal{F}_{n+1}] | \mathcal{F}_n] = \mathbb{E}[f(x)|\mathcal{F}_n] = f(x)$; i.e., $\mathbb{E}[f^{(w_{k})}(x)|\mathcal{F}_n]$ almost surely equals $f(x)$ for all $k \geq n$ and for all $x\in H$. Hence, by taking the expectation in \eqref{ineq:f} conditioned on $\mathcal{F}_n$, we have 
\begin{align}\label{f_1}
\begin{split}
&\quad \frac{2}{\lambda_{n+m}} \sum_{k=n+1}^{n+m} (1-\alpha_k)\lambda_k \left( f(x_k) - f^\star \right)\\
&\leq 
M_5 \mathbb{E} \left[ \frac{\left\| x_{n+m+1} - x_{n+1} \right\|}{\lambda_{n+m}} \bigg| \mathcal{F}_n \right]
+ \frac{\left\| x_0 - x^\star \right\|^2}{\lambda_{n+m}} \sum_{k=n+1}^{n+m} \alpha_k
 + \frac{M_1^2}{\lambda_{n+m}} \sum_{k=n+1}^{n+m} \lambda_k^2
\end{split} 
\end{align}
almost surely. Since (C5) and the monotone decreasing conditions of $(\alpha_n)_{n\in\mathbb{N}}$ and $(\lambda_n)_{n\in\mathbb{N}}$ satisfy
\begin{align}\label{stepsizes}
&\frac{\left\| x_0 - x^\star \right\|^2}{\lambda_{n+m}} \sum_{k=n+1}^{n+m} \alpha_k
\leq m \left\| x_0 - x^\star \right\|^2 \frac{\alpha_{n+1}}{\lambda_{n+m}}
\leq m(n) \left(m(n)-1 \right) \sigma \left\| x_0 - x^\star \right\|^2 \frac{\alpha_{n+1}}{\lambda_{n+1}},\nonumber\\
&\frac{M_1^2}{\lambda_{n+m}} \sum_{k=n+1}^{n+m} \lambda_k^2
\leq m M_1^2 \lambda_{n+1} \frac{\lambda_{n+1}}{\lambda_{n+m}} 
\leq m(n) \left(m(n)-1 \right) \sigma M_1^2 \lambda_{n+1},
\end{align}
Sub-assumption \ref{omega}(i), (C4), and $\lim_{n\to\infty} \lambda_n = 0$ mean that $\lim_{n\to \infty} (\left\| x_0 - x \right\|^2/\lambda_{n+m}) \sum_{k=n+1}^{n+m} \alpha_k \leq 0$ and $\lim_{n\to \infty} (M_1^2/\lambda_{n+m}) \sum_{k=n+1}^{n+m} \lambda_k^2 \leq 0$. Accordingly, Lemma \ref{lem:1} guarantees that, almost surely
\begin{align*}
\limsup_{n\to\infty} \frac{2}{\lambda_{n+m}} \sum_{k=1}^{m} (1-\alpha_{n+k})\lambda_{n+k} \left( f(x_{n+k}) - f^\star \right) \leq 0.
\end{align*}
Now, let us assume that the assertion in Lemma \ref{lem:4} does not hold; i.e., for all $\tilde{\Omega} \in \mathcal{F}$, $\mathbb{P}(\tilde{\Omega}) = 1$ and there exists $\omega \in \tilde{\Omega}$ such that $\limsup_{n\to\infty} f(x_n(\omega)) - f^\star > 0$. Accordingly, there exist $\gamma > 0$ and $n_3 \in \mathbb{N}$ such that $f(x_n(\omega)) - f^\star \geq \gamma$ for all $n \geq n_3$. The monotone decreasing conditions of $(\alpha_n)_{n\in\mathbb{N}}$ and $(\lambda_n)_{n\in\mathbb{N}}$ and $\lim_{n\to\infty}\alpha_n = 0$ thus guarantee that
\begin{align*}
0 
&\geq 
\limsup_{n\to\infty} \frac{2}{\lambda_{n+m}} \sum_{k=1}^{m} (1-\alpha_{n+k})\lambda_{n+k} \left( f(x_{n+k}(\omega)) - f^\star \right)\\
%&\geq 
%\gamma \limsup_{n\to\infty} \frac{2}{\lambda_{n+m}} \sum_{k=1}^{m} (1-\alpha_{n+k})\lambda_{n+k}\\
&\geq
\gamma \limsup_{n\to\infty} \frac{2\lambda_{n+m}}{\lambda_{n+m}} m(n) (1-\alpha_{n+1}) \geq 
2 \gamma > 0,
\end{align*}
which is a contradiction. Therefore, almost surely $\limsup_{n\to\infty} f(x_n) - f^\star \leq 0$.

Next, let us assume that Sub-assumption \ref{omega}(iv) holds. Inequality \eqref{ineq:f} thus leads to the finding that, for all $n\in \mathbb{N}$, almost surely
\begin{align*}
&\quad \frac{2(1-\alpha_{n+m})\lambda_{n+m}}{\lambda_{n+m}} \sum_{k=n+1}^{n+m} \left( f^{(w_k)}(x_k) - f^{(w_k)}(x^\star) \right)\\
&\leq M_5 \frac{\left\| x_{n+m+1} - x_{n+1} \right\|}{\lambda_{n+m}}
+ \frac{\left\| x_0 - x^\star \right\|^2 m \alpha_{n+m}}{\lambda_{n+m}}
 + \frac{M_1^2 m \lambda_{n+m}^2}{\lambda_{n+m}}.
\end{align*}
Since the definition of $\partial f^{(w_k)}$ means that $f^{(w_k)}(x_n) - f^{(w_k)}(x_k) \leq \langle x_n - x_k, \mathsf{G}^{(w_k)}(x_n) \rangle$ $(k=n+1,n+2, \ldots, n+m)$, almost surely
\begin{align}\label{ineq:5}
\begin{split}
&\quad 2(1-\alpha_{n+m}) \sum_{k=n+1}^{n+m} \left( f^{(w_k)}(x_n) - f^{(w_k)}(x^\star) \right)\\
&\leq M_5 \frac{\left\| x_{n+m+1} - x_{n+1} \right\|}{\lambda_{n+m}}
+ \frac{\left\| x_0 - x^\star \right\|^2 m \alpha_{n+m}}{\lambda_{n+m}}
+ \frac{M_1^2 m \lambda_{n+m}^2}{\lambda_{n+m}}\\
&\quad + 2M_1 (1-\alpha_{n+m}) \sum_{k=n+1}^{n+m} \left\| x_n - x_k \right\|.
\end{split}
\end{align}
Further, from $\|x_{l+1} - x_l \| \leq \|x_{l+1} - y_l\| + \|y_l - x_l\|$ and $\|x_{l+1} - y_l\|= \alpha_l \|x_0 - y_l\|$ $(l\in \mathbb{N})$, Assumption \ref{bounded} and Lemma \ref{lem:3} ensure that $\lim_{l\to\infty} \|x_{l+1} - x_l \|$ almost surely equals $0$. Hence, the triangle inequality guarantees that, for some $j\in \mathbb{N}$, $\lim_{l\to\infty} \|x_{l} - x_{l+j} \|$ almost surely equals $0$. Taking the expectation in \eqref{ineq:5} thus ensures that, for all $\epsilon > 0$, there exists $n_4 \in \mathbb{N}$ such that, for all $n \geq n_4$, almost surely
\begin{align*}
&\quad 2(1-\alpha_{n+m}) \left( f (x_n) - f^\star \right)\\
&\leq M_5 \mathbb{E} \left[ \frac{\left\| x_{n+m+1} - x_{n+1} \right\|}{m \lambda_{n+m}} \bigg| \mathcal{F}_n \right]
+ \frac{\left\| x_0 - x^\star \right\|^2 \alpha_{n+m}}{\lambda_{n+m}}
 + M_1^2 \lambda_{n+m} + 2M_1 (1-\alpha_{n+m}) \epsilon,
\end{align*}
where the left side comes from the condition that almost surely $f(x) = (1/m)\sum_{l=tm}^{(t+1)m-1} \mathbb{E}[f^{(w_l)}(x)|\mathcal{F}_{tm}] $ with $tm = n+1$ and the definition of $\mathcal{F}_n$. Hence, from Sub-assumption \ref{omega}(i), Lemma \ref{lem:1}, (C4), and $\lim_{n\to\infty} \lambda_n = \lim_{n\to\infty} \alpha_n = 0$, almost surely
\begin{align*}
2 \limsup_{n\to\infty} \left( f (x_n) - f^\star \right) \leq 2M_1 \epsilon.
\end{align*}
Therefore, the arbitrary condition of $\epsilon$ guarantees that Lemma \ref{lem:4} holds. 

Now we are in the position to prove Theorem \ref{thm:1}.

{\em Proof:}
Lemma \ref{lem:3} ensures the existence of $\bar{\Omega} \in \mathcal{F}$ such that $\mathbb{P}(\bar{\Omega}) = 1$ and $\lim_{n\to\infty}\| x_n(\omega) - T^{(i)}(x_n(\omega)) \| = 0$ for all $\omega \in \bar{\Omega}$ and for all $i\in \mathcal{I}$. Moreover, Lemma \ref{lem:4} means that there exists $\hat{\Omega}\in \mathcal{F}$ such that $\mathbb{P}(\hat{\Omega}) = 1$ and $\limsup_{n\to\infty} f(x_n(\omega)) \leq f^\star$ for all $\omega \in \hat{\Omega}$. Now, let $\omega \in \bar{\Omega} \cap \hat{\Omega}$ and let $x^* \in \mathcal{W}(x_n(\omega))_{n\in \mathbb{N}}$. Assumption \ref{bounded} and $\mathbb{P}(\bar{\Omega} \cap \hat{\Omega}) = 1$ guarantee the existence of a weak sequential cluster point of $(x_n(\omega))_{n\in\mathbb{N}}$. Then, there exists $(x_{n_l}(\omega))_{l\in \mathbb{N}} \subset (x_n(\omega))_{n\in \mathbb{N}}$ such that it converges weakly to $x^* \in H$. Here, let us fix $i\in \mathcal{I}$ arbitrarily and assume that $x^* \notin \mathrm{Fix}(T^{(i)})$. From Opial's lemma \cite[Lemma 3.1]{opial}, 
\begin{align*}
\liminf_{l\to\infty} \left\| x_{n_l}(\omega) - x^* \right\|
< \liminf_{l\to\infty} \left\| x_{n_l}(\omega) - T^{(i)} (x^*) \right\|,
\end{align*}
which, together with $\omega \in \bar{\Omega}$ and (A1), means that 
\begin{align*}
\liminf_{l\to\infty} \left\| x_{n_l}(\omega) - x^* \right\|
< \liminf_{l\to\infty} \left\| T^{(i)}(x_{n_l}(\omega)) - T^{(i)}(x^*) \right\|
\leq \liminf_{l\to\infty} \left\| x_{n_l}(\omega)- x^* \right\|.
\end{align*}
This is a contradiction. Therefore, $x^* \in \mathrm{Fix}(T^{(i)})$ for all $i\in \mathcal{I}$; i.e., $x^* \in X$. Furthermore, the weakly lower semicontinuity of $f$ \cite[Theorem 9.1]{b-c} leads to the finding that
\begin{align*}
f(x^*) \leq \liminf_{l\to \infty} f\left(x_{n_l}(\omega)\right)
\leq \limsup_{n\to \infty} f\left(x_{n}(\omega)\right) \leq f^\star.
\end{align*}
That is, $x^* \in X^\star$. This completes the proof.

\subsection{Convergence rate analysis of Algorithm \ref{algo:1}} 
The following proposition establishes the rate of convergence for Algorithm \ref{algo:1}.

\begin{prop}\label{thm:2}
Suppose that the assumptions in Theorem \ref{thm:1} hold and that $(x_n)_{n\in \mathbb{N}}$ is the sequence generated by Algorithm \ref{algo:1}. Then, there exist $N_i \in \mathbb{R}$ ($i=1,2$) such that, for all $i\in \mathcal{I}$ and for all $n\in \mathbb{N}$, almost surely
\begin{align*}%\label{rate:1_1}
\left\| x_n - T^{(i)} (x_n) \right\|
\leq 
\sqrt{N_1 \alpha_{n} + N_2 \lambda_{n}}. 
\end{align*}
Moreover, under Sub-assumption \ref{omega}(iii), if there exists $k_0 \in \mathbb{N}$ such that $f (x_n) \geq f^\star$ almost surely for all $n \geq k_0$, then there exist $k_1 \in \mathbb{N}$ and $N_i \in \mathbb{R}$ ($i=3,4,5$) such that, for all $n \geq \max\{k_0, k_1\}$, almost surely
\begin{align}\label{rate:1_2}
\frac{1}{m} \sum_{k=n+1}^{n+m} f(x_k) - f^\star
\leq 
N_3 \frac{o(\lambda_{n+m})}{\lambda_{n+m}} + N_4 \lambda_{n} + N_5 \frac{\alpha_n}{\lambda_{n}}.
\end{align}
Under Sub-assumption \ref{omega}(iv), there exist $k_2 \in \mathbb{N}$ and $N_i \in \mathbb{R}$ ($i=6,7,8,9,10$) such that, for all $n \geq k_2$, almost surely
\begin{align}\label{rate:1_3}
f(x_n) - f^\star 
&\leq 
N_6 \frac{o(\lambda_{n+m})}{\lambda_{n+m}} 
+ N_7 \lambda_n + N_8 \frac{\alpha_{n}}{\lambda_{n}}
+ \sqrt{N_9 \alpha_{n} + N_{10} \lambda_{n}}.
\end{align}
\end{prop}

Here, let us compare the stochastic first-order method with random constraint projection \cite{bert2015} with Algorithm \ref{algo:1}. In \cite{bert2015}, the problem
\begin{align}\label{problem}
\text{minimize } f(x) := \mathbb{E}\left[ f^{(v)}(x) \right]
\text{ subject to } x \in C := \bigcap_{i=1}^M C^{(i)}
\end{align}
was discussed \cite[(1)--(3)]{bert2015}, where $f^{(v)} \colon \mathbb{R}^N \to \mathbb{R}$ is a convex function of $x$ involving a random variable $v$, and $C^{(i)} \subset \mathbb{R}^N$ $(i=1,2,\ldots,M)$ is a nonempty, closed convex set onto which the metric projection $P^{(i)}$ can be efficiently computed. The following stochastic first-order method \cite[Algorithm 1, (9)]{bert2015} was presented for solving problem \eqref{problem}: given $x_k \in\mathbb{R}^N$,
\begin{align}\label{bert}
\begin{split}
&z_k := x_k - \alpha_k \mathsf{G}^{(v_k)}(\bar{x}_k),\\
&x_{k+1} := z_k - \beta_k \left( z_k - \mathsf{P}^{(w_k)} (z_k) \right),\text{ with } \bar{x}_k = x_k \text{ or } \bar{x}_k = x_{k+1},
\end{split}
\end{align}
where $\mathsf{P}^{(w)}$ stands for the stochastic metric projection onto $C^{(w)}$, and $(\alpha_k)_{k\in \mathbb{N}}, (\beta_k)_{k\in \mathbb{N}} \subset (0,\infty)$. Under certain assumptions, Algorithm \eqref{bert} converges almost surely to a random point in the solution set of problem \eqref{problem} \cite[Theorem 1]{bert2015}. Theorem 2 in \cite{bert2015} implies that, under certain assumptions, Algorithm \eqref{bert} with $\alpha_k = 1/\sqrt{k}$ and $\beta_k := \beta > 0$ $(k\in \mathbb{N})$ satisfies 
\begin{align*}
\mathbb{E}\left[ f \left( \frac{1}{k} \sum_{t=1}^k P_C (x_t) \right) \right]
= f^* + O\left( \frac{1}{\sqrt{k}} \right),
\text{ } 
\mathbb{E}\left[
\mathrm{d}\left( \frac{1}{k} \sum_{t=1}^k x_t, C \right)^2
\right] 
= O\left( \frac{\log (k+1)}{k} \right),
\end{align*}
where $f^*$ is the optimal value of problem \eqref{problem} and $\mathrm{d}(x,C) := \inf_{y\in C} \|x-y\|$ $(x\in \mathbb{R}^N)$. 

Meanwhile, Algorithm \ref{algo:1} can be applied to problem \eqref{problem} even when $C^{(i)}$ is not always simple in the sense that $P^{(i)}$ cannot be easily computed (see Section \ref{sec:5} for an example of problem \eqref{problem} when $C^{(i)}$ is not simple). Theorem \ref{thm:1} guarantees that any weak sequential cluster point of $(x_n)_{n\in \mathbb{N}}$ generated by Algorithm \ref{algo:1} almost surely belongs to the solution set of Problem \ref{prob:1} including problem \eqref{problem}. Proposition \ref{thm:2} implies that Algorithm \ref{algo:1} with $\lambda_n := 1/n^a$ and $\alpha_n := 1/n^b$ $(n\geq 1)$, where $a \in (0,1/2)$ and $b\in (a,1-a)$, satisfies, for all $i\in \mathcal{I}$,
\begin{align}\label{rate:1_1}
\left\| x_n - T^{(i)} (x_n) \right\| 
= O \left(\frac{1}{\sqrt{n^{a}}} \right).
\end{align}
Moreover, \eqref{rate:1_2} implies that, under the assumptions in Proposition \ref{thm:2} and the condition $o(\lambda_{n}) = 1/n^c$, where $c > a$, 
\begin{align}\label{Rate}
\frac{1}{m} \sum_{k=n+1}^{n+m} f(x_k) - f^\star
=
O \left( \frac{1}{n^{\min\{ a,b-a,c-a \}}} \right),
\end{align}
while \eqref{rate:1_3} implies that
\begin{align}\label{Rate:1} 
f(x_n) - f^\star 
= 
O \left( \frac{1}{n^{\min\{a/2, b-a, c-a\}}} \right).
\end{align}

{\em Proof:}
From \eqref{xnyn2}, the monotone decreasing conditions of $(\alpha_n)_{n\in\mathbb{N}}$ and $(\lambda_n)_{n\in\mathbb{N}}$ with $\lim_{n\to\infty}\alpha_n = 0$, and the almost sure boundedness of $(x_n)_{n\in\mathbb{N}}$, there exist $N_i \in \mathbb{R}$ $(i=1,2)$ such that, for all $n\in \mathbb{N}$, almost surely $\mathbb{E} [\|x_{n} - y_{n}\|^2| \mathcal{F}_{n} ] \leq N_1 \alpha_{n} + N_2 \lambda_{n}$. Accordingly, \eqref{xnt} ensures that, for all $n\in \mathbb{N}$, almost surely
\begin{align*}
\mathbb{E} \left[ \left\| x_{n} - \mathsf{T}^{(w_n)} (x_{n}) \right\|^2 \bigg|
\mathcal{F}_n \right]
&\leq 
2 \left(N_1 \alpha_{n} + N_2 \lambda_{n} \right)
+ 2 M_1^2 \lambda_{n}^2.
\end{align*}
Sub-assumption \ref{omega}(ii) guarantees the existence of $N_3 \in \mathbb{R}$ such that, for $i\in \mathcal{I}$ and for all $n\in \mathbb{N}$, $\| x_n - T^{(i)}(x_n) \|^2 \leq N_3 \mathbb{E}[\| x_n - \mathsf{T}^{(w_n)}(x_n) \|^2 |\mathcal{F}_n]$ holds almost surely. This means that, for all $i\in \mathcal{I}$ and for all $n\in \mathbb{N}$, almost surely
\begin{align*}
\left\| x_{n} - T^{(i)} (x_{n}) \right\|^2 
\leq 
2N_3 \left(N_1 \alpha_{n} + N_2 \lambda_{n} \right)
+ 2 N_3 M_1^2 \lambda_{n}.
\end{align*}
Lemma \ref{lem:1} and the monotone decreasing condition of $(\lambda_n)_{n\in\mathbb{N}}$ with $\lim_{n\to\infty} \lambda_n = 0$ guarantee that, for all $n\in \mathbb{N}$, almost surely
\begin{align}\label{key}
\mathbb{E} \left[\frac{\left\| x_{n+m+1} - x_{n+1} \right\|}{\lambda_{n+m}} \Bigg| \mathcal{F}_{n} \right]
= \frac{\mathbb{E} \left[ \left\| x_{n+m+1} - x_{n+1} \right\| | \mathcal{F}_{n} \right]}{\lambda_{n+m}}
= \frac{o(\lambda_{n+m})}{\lambda_{n+m}}.
\end{align}

Let us assume that Sub-assumption \ref{omega}(iii) holds. From \eqref{f_1} and \eqref{key}, for all $n \in \mathbb{N}$, almost surely 
\begin{align*}
\frac{2}{\lambda_{n+m}} \sum_{k=n+1}^{n+m} (1-\alpha_k)\lambda_k \left( f(x_k) - f^\star \right)
\leq 
M_5 \frac{o(\lambda_{n+m})}{\lambda_{n+m}}
+ \frac{\left\| x_0 - x^\star \right\|^2}{\lambda_{n+m}} \sum_{k=n+1}^{n+m} \alpha_k
 + \frac{M_1^2}{\lambda_{n+m}} \sum_{k=n+1}^{n+m} \lambda_k^2.
\end{align*}
Hence, \eqref{stepsizes} guarantees that there exist $N_4, N_5 \in \mathbb{R}$ such that, for all $n \in \mathbb{N}$, almost surely
\begin{align*}
\frac{2}{\lambda_{n+m}} \sum_{k=n+1}^{n+m} (1-\alpha_k)\lambda_k \left( f(x_k) - f^\star \right)
\leq 
M_5 \frac{o(\lambda_{n+m})}{\lambda_{n+m}} 
+ N_4 \frac{\alpha_{n}}{\lambda_{n}}
 + N_5 \lambda_{n}.
\end{align*}
Since $(\alpha_n)_{n\in \mathbb{N}}$ converges to $0$, there exists $n_5 \in \mathbb{N}$ such that $1 - \alpha_n \geq 1/2$ for all $n \geq n_5$. From the monotone decreasing condition of $(\lambda_n)_{n\in\mathbb{N}}$ and the existence of $n_6 \in \mathbb{N}$ such that almost surely $f(x_n) - f^\star \geq 0$ for all $n \geq n_6$, we have for all $n \geq k_0 := \max \{n_5, n_6\}$, 
\begin{align*}
\frac{1}{m} \sum_{k=n+1}^{n+m} f(x_k) - f^\star 
\leq 
\frac{M_5}{m} \frac{o(\lambda_{n+m})}{\lambda_{n+m}}
+ \frac{N_4}{m} \frac{\alpha_{n}}{\lambda_{n}} + \frac{N_5}{m} \lambda_n
\end{align*}
almost surely. 

Let us assume that Sub-assumption \ref{omega}(iv) holds. Then, \eqref{ineq:5} guarantees that, for all $n \in \mathbb{N}$, almost surely
\begin{align*}
&\quad 2(1-\alpha_{n+m}) \sum_{k=n+1}^{n+m} \left( f^{(w_k)}(x_n) - f^{(w_k)}(x^\star) \right)\\
&\leq M_5 \frac{\left\| x_{n+m+1} - x_{n+1} \right\|}{\lambda_{n+m}}
+ N_4 \frac{\alpha_{n}}{\lambda_{n}}
+ N_5 \lambda_{n}
+ 2M_1 (1-\alpha_{n+m}) \sum_{k=n+1}^{n+m} \left\| x_n - x_k \right\|.
\end{align*}
From \eqref{xnyn3}, $\| x_{n+1} - y_n \| = \alpha_n \|x_0 - y_n \|$ $(n\in \mathbb{N})$, and the triangle inequality, there exist $\bar{n}_1 \in \mathbb{N}$ and $\bar{N}_i \in \mathbb{R}$ $(i=1,2,3)$ such that, for all $n \geq \bar{n}_1$, almost surely
\begin{align*}
\left\| x_{n} - x_{n+1} \right\| 
&\leq \left\| x_{n} - y_n \right\| +\left\| y_n - x_{n+1} \right\|\\
&\leq \sqrt{\bar{N}_1 \alpha_n + \bar{N}_2 \lambda_n} + \bar{N}_3 \alpha_n, 
\end{align*}
which, together with the triangle inequality and the monotone decreasing conditions of $(\alpha_n)_{n\in\mathbb{N}}$ and $(\lambda_n)_{n\in\mathbb{N}}$, means that, for all $n \geq \bar{n}_1$, almost surely
\begin{align*}
\sum_{k=n+1}^{n+m} \left\| x_n - x_k \right\|
&\leq \sum_{j=0}^{m-1} (m- j) \|x_{n+j} - x_{n+j+1} \|\\
&\leq \sum_{j=0}^{m-1} (m- j) \left(\sqrt{\bar{N}_1 \alpha_{n+j} + \bar{N}_2 \lambda_{n+j}} + \bar{N}_3 \alpha_{n+j} \right)\\
&\leq \frac{m(m+1)}{2} \left(\sqrt{\bar{N}_1 \alpha_{n} + \bar{N}_2 \lambda_{n}} + \bar{N}_3 \alpha_{n} \right).
\end{align*}
Since $(\alpha_n)_{n\in\mathbb{N}}$ converges to $0$, there exist $a \in (0,1)$ and $\bar{n}_2 \in \mathbb{N}$ such that $a \leq 2(1-\alpha_n)$ for all $n \geq \bar{n}_2$. From \eqref{key}, for all $n > \max\{ n_2, \bar{n}_1, \bar{n}_2 \}$, almost surely
\begin{align*}
m\left(f(x_n) - f^\star \right) 
&\leq \frac{M_5}{a} \frac{o(\lambda_{n+m})}{\lambda_{n+m}}
+ \frac{N_4}{a} \frac{\alpha_{n}}{\lambda_{n}}
+ \frac{N_5}{a} \lambda_{n}\\
&\quad + \frac{m(m+1)M_1}{2} \left(\sqrt{\bar{N}_1 \alpha_{n} + \bar{N}_2 \lambda_{n}} + \bar{N}_3 \alpha_{n} \right),
\end{align*}
which, together with (C4) (i.e., there exists $M \in \mathbb{R}$ such that $\alpha_n \leq M \lambda_n$ for all $n\in\mathbb{N}$), completes the proof.

The following remark is made regarding Proposition \ref{thm:3}.

\begin{rem}\label{rem:1}
{\em From a discussion similar to the ones for obtaining \eqref{ineq:1} and \eqref{ineq:2}, there exist $\bar{M}_i \in \mathbb{R}$ $(i=1,2)$ such that, for all $n > n_2$, almost surely
\begin{align}\label{key_1}
\mathbb{E} \left[\frac{\left\| x_{n+m(n)+1} - x_{n+1} \right\|}{\lambda_{n+m(n)}} \Bigg| \mathcal{F}_{n} \right]
\leq
\bar{M}_1 \prod_{k=n_2}^{n} \left(1 - \alpha_{k+m(k)} \right) + \bar{M}_2 N(n),
\end{align}
where $N(n):= \max \{ (1/\alpha_{k+m(k)}) | (1/\lambda_k) - (1/\lambda_{k+m(k)}) |, (1/\alpha_{k+m(k)}) | (1/\lambda_{k+m(k)}) - (1/\lambda_{k+m(k)-1}) |, (1 /\lambda_{k+m(k)}) | 1 - \alpha_k/\alpha_{k+m(k)}| \colon k=n,n-1,\ldots,n_2 \}$. Accordingly, \eqref{key} can be replaced with \eqref{key_1}.}
\end{rem}

\section{Stochastic proximal point algorithm for nonsmooth convex optimization}
\label{sec:4}
This section presents the convergence analysis of the following proximal-type algorithm for solving Problem \ref{prob:1}.

\begin{algorithm} 
\caption{Stochastic proximal point algorithm for solving Problem \ref{prob:1}} 
\label{algo:2} 
\begin{algorithmic}[1] 
\REQUIRE
$n\in \mathbb{N}$, $(\alpha_n)_{n\in\mathbb{N}}, (\gamma_n)_{n\in \mathbb{N}} \subset (0,\infty)$.
\STATE
$n \gets 0$, 
$x_0 \in H$
\LOOP 
 \STATE 
 $y_{n} := \mathsf{T}^{(w_n)} \left(\mathsf{Prox}_{\gamma_n f^{(w_n)}}(x_n) \right)$
 \STATE 
 $x_{n+1} := \alpha_n x_0 + (1-\alpha_n) y_n$ 
 \STATE
 $n \gets n+1$
\ENDLOOP 
\end{algorithmic}
\end{algorithm}

Algorithms \ref{algo:1} and \ref{algo:2} are based on the Halpern fixed point algorithm \cite{halpern,wit}. In contrast to Algorithm \ref{algo:1}, Algorithm \ref{algo:2} uses the approach of proximal point algorithms \cite[Chapter 27]{b-c}, \cite{bacak2014,bert2011,lions,martinet1970,rock1976,bert2015} that optimize nonsmooth, convex functions over the whole space.

\subsection{Assumptions for convergence analysis of Algorithm \ref{algo:2}}
\label{subsec:4.1}
Let us consider Problem \ref{prob:1} under (A1), (A2), and (A4) defined as follows.
\begin{enumerate}
\item[(A4)] $\mathrm{Prox}_{\gamma f^{(i)}}$ $(\gamma > 0, i\in \mathcal{I})$ can be efficiently computed.
\end{enumerate} 
Tables 10.1 and 10.2 in \cite{comb2011} provide several examples of convex functions for which proximity operators can be computed within a finite number of arithmetic operations.

The conditions of the step-size sequences in Algorithm \ref{algo:2} are as follows.

\begin{assum}\label{stepsize2}
Let $\sigma \geq 1$. The step-size sequences $(\alpha_n)_{n\in \mathbb{N}} \subset (0,1)$ and $(\gamma_n)_{n\in \mathbb{N}} \subset (0,1)$, which are monotone decreasing and converge to $0$, satisfy the following conditions:
\begin{align*}
&\text{{\em (C1)}} \sum_{n=0}^\infty \alpha_n = \infty, \text{ }
\text{{\em (C2)}} \lim_{n\to\infty} \frac{1}{\alpha_{n+1}} \left| \frac{1}{\gamma_{n+1}} - \frac{1}{\gamma_n} \right| = 0, \text{ } 
\text{{\em (C3)}} \lim_{n\to\infty} \frac{1}{\gamma_{n+1}} \left| 1 - \frac{\alpha_n}{\alpha_{n+1}} \right| = 0,\\
&\text{{\em (C4)}} \lim_{n\to\infty} \frac{\alpha_n}{\gamma_n} = 0, \text{ }
\text{{\em (C5)}} \lim_{n\to\infty} \frac{1}{\alpha_{n+1}}\frac{|\gamma_{n+1} - \gamma_n|}{\gamma_{n+1}^2} = 0, \text{ }
\text{{\em (C6)}} \frac{\alpha_n}{\alpha_{n+1}}, \frac{\lambda_n}{\lambda_{n+1}} \leq \sigma
\text{ } (n\in \mathbb{N}). 
%\text{ } \lim_{n\to\infty} \frac{\lambda_n}{\lambda_{n+1}} = 1.
\end{align*}
\end{assum}
Examples of $(\alpha_n)_{n\in \mathbb{N}}$ and $(\gamma_n)_{n\in \mathbb{N}}$ satisfying Assumption \ref{stepsize2} are $\gamma_n := 1/(n+1)^a$ and $\alpha_n := 1/(n+1)^b$ $(n\in \mathbb{N})$, where $a\in (0,1/2)$, $b\in (a,1-a)$, and $a+b < 1$.

The convergence of Algorithm \ref{algo:2} depends on the following assumption.
\begin{assum}\label{bounded2}
The sequence $(w_n)_{n\in \mathbb{N}}$ satisfies Assumption \ref{omega}, where $\lambda_n$ is replaced with $\gamma_n$. The sequence $(y_n)_{n\in \mathbb{N}}$ is almost surely bounded.
\end{assum}

A similar discussion to the one for defining \eqref{y_n} implies that, if there exists a simple, bounded, closed convex set $C \supset X$, then $y_n$ $(n\in \mathbb{N})$ in Algorithm \ref{algo:2} can be replaced with
\begin{align}\label{y_n2}
y_n := P_C \left[\mathsf{T}^{(w_n)} \left( \mathsf{Prox}_{\gamma_n f^{(w_n)}}(x_n) \right) \right],
\end{align}
which implies the boundedness of $(y_n)_{n\in \mathbb{N}}$.

\subsection{Convergence analysis of Algorithm \ref{algo:2}}
\label{subsec:4.2}
\text{}

\begin{thm}\label{thm:3}
Suppose that Assumptions (A1), (A2), (A4), \ref{stepsize2}, and \ref{bounded2} hold, and let $(x_n)_{n\in \mathbb{N}}$ be the sequence generated by Algorithm \ref{algo:2}. Then, any weak sequential cluster point of $(x_n)_{n\in\mathbb{N}}$ almost surely belongs to the solution set of Problem \ref{prob:1}.
\end{thm}

The proof starts with the following lemma.

\begin{lem}\label{lem:1_1}
Suppose that the assumptions in Theorem \ref{thm:3} hold. Then, almost surely
\begin{align*}
\lim_{n\to\infty} \mathbb{E} \left[\frac{\| x_{n+m+1} - x_{n+1} \|}{\gamma_{n+m}} 
\bigg| \mathcal{F}_n \right] = 0.
\end{align*} 
\end{lem}

{\em Proof:}
Sub-assumption \ref{omega}(i) ensures that, for all $n\in \mathbb{N}$, there exists $m(n) \in \mathbb{N}$ such that $\limsup_{n\to \infty} m(n) < \infty$, $\mathsf{T}^{(w_{n+m})} = \mathsf{T}^{(w_n)}$, and $f^{(w_{n+m})} = f^{(w_n)}$ almost surely. Accordingly, (A1) and the triangle inequality ensure that, for all $n\geq n_0$, almost surely
\begin{align*}
\left\| y_{n+m} - y_n \right\|
%= &\left\| T^{(w_n)} \left(\mathrm{Prox}_{\gamma_{n+m} f^{(w_{n})}}
% (x_{n+m}) \right) 
% - T^{(w_n)} \left(\mathrm{Prox}_{\gamma_{n} f^{(w_{n})}}
% (x_{n}) \right) \right\|\\
%\leq & \left\| \mathrm{Prox}_{\gamma_{n+m} f^{(w_{n})}}
% (x_{n+m}) 
% - \mathrm{Prox}_{\gamma_{n} f^{(w_{n})}}
% (x_{n})\right\|\\
\leq & \left\| \mathsf{Prox}_{\gamma_{n+m} f^{(w_{n})}}
 (x_{n+m}) 
 - \mathsf{Prox}_{\gamma_{n+m} f^{(w_{n})}}
 (x_{n})\right\|\\
 &+ \left\| \mathsf{Prox}_{\gamma_{n+m} f^{(w_{n})}}
 (x_{n}) 
 - \mathsf{Prox}_{\gamma_{n} f^{(w_{n})}}
 (x_{n})\right\|, 
\end{align*}
which, together with Proposition \ref{prop:1}(ii), means that
\begin{align*}
\left\| y_{n+m} - y_n \right\|
\leq \left\| x_{n+m} - x_n \right\| + \left\| \mathsf{Prox}_{\gamma_{n+m} f^{(w_{n})}}(x_{n}) 
- \mathsf{Prox}_{\gamma_{n} f^{(w_{n})}}(x_{n})\right\|.
\end{align*}
Put $z_n := \mathsf{Prox}_{\gamma_{n} f^{(w_{n})}}(x_{n})$ and $\bar{z}_n := \mathsf{Prox}_{\gamma_{n+m} f^{(w_{n})}}(x_{n})$ $(n\in\mathbb{N})$. Proposition \ref{prop:1}(i) thus means that $(x_{n} - z_n)/\gamma_{n} \in \partial f^{(w_{n})}(z_n)$ and $(x_n - \bar{z}_n)/\gamma_{n+m} \in \partial f^{(w_{n})}(\bar{z}_n)$ $(n\in\mathbb{N})$. Hence, the monotonicity of $\partial f^{(w_n)}$ implies that, for all $n\in \mathbb{N}$, $\langle z_n - \bar{z}_n, (x_n - z_n)/\gamma_n - (x_n - \bar{z}_n)/\gamma_{n+m} \rangle \geq 0$, which means that 
\begin{align*}
&\frac{1}{\gamma_n \gamma_{n+m}}
\big\{
\left\langle z_n - \bar{z}_n, (\gamma_{n+m} - \gamma_n) x_n \right\rangle
+
\left\langle z_n - \bar{z}_n, -\gamma_{n+m} (z_n - \bar{z}_n) \right\rangle\\
&\quad 
+
\left\langle z_n - \bar{z}_n, (\gamma_n-\gamma_{n+m}) \bar{z}_n \right\rangle
\big\} \geq 0.
\end{align*}
Accordingly, for all $n\in \mathbb{N}$, $\left\| z_n - \bar{z}_n \right\| \leq (|\gamma_{n+m} - \gamma_n|/\gamma_{n+m}) (\|x_n\| + \|\bar{z}_n\|)$. Thus, for all $n\in \mathbb{N}$, almost surely
\begin{align*}
\left\| y_{n+m} - y_n \right\| 
\leq \left\| x_{n+m} - x_n \right\| + \frac{|\gamma_{n+m} - \gamma_n|}{\gamma_{n+m}}
\left( \|x_n\| + \|\bar{z}_n\| \right).
\end{align*}
A discussion similar to the one for obtaining \eqref{xn} guarantees that, for all $n\in \mathbb{N}$, almost surely
\begin{align*}
\left\| x_{n+m+1} - x_{n+1} \right\| 
%&\leq \left(1 - \alpha_{n+m} \right) \left\{ \left\| x_{n+m} - x_n \right\| + \frac{|\gamma_{n+m} - \gamma_n|}{\gamma_{n+m}}
%\left( \|x_n\| + \|\bar{z}_n\| \right) \right\}\\ 
%&\quad + \left| \alpha_{n+m} - \alpha_n \right| \left\|x_0 - y_n\right\|\\
&\leq \left(1 - \alpha_{n+m} \right)\left\| x_{n+m} - x_n \right\| + 
 \frac{|\gamma_{n+m} - \gamma_n|}{\gamma_{n+m}}
\left( \|x_n\| + \|\bar{z}_n\| \right)\\
&\quad + \left| \alpha_{n+m} - \alpha_n \right| \left\|x_0 - y_n\right\|.
\end{align*}
Therefore, the same discussion as for \eqref{ineq:1} implies that, for all $n\in \mathbb{N}$, almost surely
\begin{align*}%\label{ineq:1}
\frac{\left\| x_{n+m+1} - x_{n+1} \right\|}{\gamma_{n+m}}
\leq&
\left(1 - \alpha_{n+m} \right) \frac{\left\| x_{n+m} - x_n \right\|}{\gamma_{n+m-1}}
+ \alpha_{n+m} \frac{1}{\gamma_{n+m}}\left| 1 - \frac{\alpha_n}{\alpha_{n+m}} \right| \left\|x_0 - y_n\right\|\\ 
 &+ \alpha_{n+m} \frac{1}{\alpha_{n+m}} \left| \frac{1}{\gamma_{n+m}} - \frac{1}{\gamma_{n+m-1}} \right| \left\| x_{n+m} - x_n \right\|\\
 &+ \alpha_{n+m} \frac{1}{\alpha_{n+m}} \frac{|\gamma_{n+m} - \gamma_n|}{\gamma_{n+m}^2}\left( \|x_n\| + \|\bar{z}_n\| \right).
\end{align*}
Therefore, the proof of Lemma \ref{lem:1}, Sub-assumption \ref{omega}(i), and Assumptions \ref{stepsize2} and \ref{bounded2} lead to the assertion in Lemma \ref{lem:1_1}. This completes the proof.

\begin{lem}\label{lem:2_1}
Suppose that the assumptions in Theorem \ref{thm:3} hold and $z_n := \mathsf{Prox}_{\gamma_n f^{(w_n)}}(x_n)$ for all $n\in \mathbb{N}$. Then, almost surely
\begin{align*}
\lim_{n\to\infty} \mathbb{E} \left[\left\|x_{n} - z_{n} \right\|^2
\Big| \mathcal{F}_n \right] = 0 \text{ and }
\lim_{n\to\infty} \mathbb{E}\left[ \left\| x_{n} - \mathsf{T}^{(w_n)} (x_{n}) \right\|^2 \bigg| \mathcal{F}_n \right] = 0.
\end{align*}
\end{lem}

{\em Proof:}
Choose $x\in X$ and $n\in \mathbb{N}$ arbitrarily and define $z_k := \mathsf{Prox}_{\gamma_k f^{(w_k)}}(x_k)$ $(k\in \mathbb{N})$. Proposition \ref{prop:1}(i) thus ensures that, for all $k\in \mathbb{N}$, $\langle x - z_k, x_k - z_k \rangle \leq \gamma_k ( f^{(w_k)} (x) - f^{(w_k)} (z_k) )$, which, together with $\langle x,y \rangle = (1/2)(\|x\|^2 + \|y\|^2 - \|x-y\|^2)$ $(x,y\in H)$, means that
\begin{align*}
\left\| z_k - x \right\|^2 
\leq \left\| x_k - x \right\|^2 - \left\| z_k - x_k \right\|^2
+ 2 \gamma_k \left( f^{(w_k)} (x) - f^{(w_k)} (z_k) \right).
\end{align*}
Since the convexity of $\|\cdot\|^2$ and (A1) mean that, for all $k\in \mathbb{N}$, $\| x_{k+1} - x \|^2 \leq \alpha_k \| x_{0} - x \|^2 + \|z_k - x \|^2 - (1-\alpha_k) \| z_k - \mathsf{T}^{(w_k)}(z_k) \|^2$, we also have, for all $k\in \mathbb{N}$,
\begin{align}
\left\| x_{k+1} - x \right\|^2
&\leq 
\alpha_k \left\| x_{0} - x \right\|^2 + \left\| x_k - x \right\|^2 - \left\| z_k - x_k \right\|^2
+ 2 \gamma_k \left( f^{(w_k)} (x) - f^{(w_k)} (z_k) \right)\nonumber\\
&\quad - (1-\alpha_k) \left\| z_k - \mathsf{T}^{(w_k)}(z_k) \right\|^2.\label{keyi}
\end{align}
Furthermore, the definition of $\partial f^{(i)}$ $(i\in \mathcal{I})$ and Proposition \ref{prop:1}(iii) imply that there exists $K_1 \in \mathbb{R}$ such that
\begin{align*}
\left\| x_{k+1} - x \right\|^2
&\leq 
\alpha_k \left\| x_{0} - x \right\|^2 + \left\| x_k - x \right\|^2 - \left\| z_k - x_k \right\|^2
+ 2 K_1 \gamma_k \left\| x - z_k \right\|\\
&\quad - (1-\alpha_k) \left\| z_k - \mathsf{T}^{(w_k)}(z_k) \right\|^2.
\end{align*}
Accordingly, 
\begin{align}
&\left\| z_{n+1} - x_{n+1} \right\|^2
\leq 
\left\| x_{0} - x \right\|^2 \sum_{k=n+1}^{n+m} \alpha_k
+ 2 K_1 \sum_{k=n+1}^{n+m} \gamma_k \left\| x - z_k \right\|\nonumber\\
&\quad\quad\quad\quad\quad\quad\quad\quad 
 + \gamma_{n+m} \left( \left\| x_{n+1} - x \right\| + \left\| x_{n+m+1} - x \right\| \right) 
\frac{\left\| x_{n+1} - x_{n+m+1} \right\|}{\gamma_{n+m}},\label{constant1}\\
&(1 - \alpha_{n+1}) \left\| z_{n+1} - \mathsf{T}^{(w_{n+1})}(z_{n+1}) \right\|^2
\leq 
\left\| x_{0} - x \right\|^2 \sum_{k=n+1}^{n+m} \alpha_k
+ 2 K_1 \sum_{k=n+1}^{n+m} \gamma_k \left\| x - z_k \right\|\nonumber\\
&\quad\quad\quad\quad\quad\quad\quad\quad 
+ \gamma_{n+m} \left( \left\| x_{n+1} - x \right\| + \left\| x_{n+m+1} - x \right\| \right) \frac{\left\| x_{n+1} - x_{n+m+1} \right\|}
{\gamma_{n+m}}.\label{constant2}
\end{align}
Therefore, from a discussion similar to the one for obtaining \eqref{xnyn}, Assumptions \ref{stepsize2} and \ref{bounded2}, and Lemma \ref{lem:1_1} lead to $\lim_{n\to\infty} \mathbb{E} [ \| z_{n} - x_{n} \|^2 | \mathcal{F}_n] = 0$ almost surely and $\lim_{n\to\infty} \mathbb{E} [ \| z_{n} - \mathsf{T}^{(w_n)}(z_{n}) \|^2 | \mathcal{F}_n ] = 0$ almost surely. Since (A1) and $\|x+y\|^2 \leq 2\|x\|^2 + 2 \|y\|^2$ $(x,y\in H)$ guarantee that
\begin{align}\label{constant3}
\begin{split}
\left\| x_n - \mathsf{T}^{(w_n)} (x_n) \right\|^2
&\leq 
2 \left\| x_n - z_n \right\|^2 + 2 \left\| z_n - \mathsf{T}^{(w_n)} (x_n) \right\|^2\\
&\leq
6 \left\| x_n - z_n \right\|^2 
+ 4 \left\| z_n - \mathsf{T}^{(w_n)} (z_n) \right\|^2,
\end{split}
\end{align}
we have $\lim_{n\to\infty} \mathbb{E}[\| x_{n} - \mathsf{T}^{(w_n)}(x_{n})\|^2| \mathcal{F}_n] = 0$ almost surely. This completes the proof.

Lemma \ref{lem:2_1} leads to the following.
 
\begin{lem}\label{lem:3_1}
Suppose that the assumptions in Theorem \ref{thm:3} hold. Then, for all $i\in \mathcal{I}$, almost surely
\begin{align*}
\lim_{n\to\infty} \left\| x_n - T^{(i)}(x_n)\right\| = 0
\text{ and }
\lim_{n\to\infty} \left\| x_n - z_n \right\| = 0.
\end{align*}
\end{lem}

{\em Proof:}
The same discussion as for proving Lemma \ref{lem:3} guarantees that $\lim_{n\to\infty} \| x_n - T^{(i)}(x_n) \| = 0$ $(i\in \mathcal{I})$ almost surely. Lemma \ref{lem:1_1} ensures that $(\|x_{n+m+1} - x_{n+1}\|/\gamma_{n+m})_{n\in \mathbb{N}}$ almost surely is bounded. Hence, \eqref{constant1} and $\lim_{n\to\infty} \alpha_n = \lim_{n\to\infty}\gamma_n = 0$ guarantee that $\lim_{n\to\infty}\|x_n - z_n\|$ almost surely equals $0$. This completes the proof.

Lemma \ref{lem:3_1} leads to the following.
\begin{lem}\label{lem:4_1}
Suppose that the assumptions in Theorem \ref{thm:3} hold. Then, almost surely
\begin{align*}
\limsup_{n\to\infty} f(x_n) \leq f^\star := \min_{x\in X}f(x).
\end{align*}
Moreover, any weak sequential cluster point of $(x_n)_{n\in\mathbb{N}}$ almost surely belongs to $X^\star := \{ x^\star \in X \colon f(x^\star) = f^\star \}$.
\end{lem}

{\em Proof:}
Choose $x^\star \in X^\star$ and $n\in \mathbb{N}$ arbitrarily. Inequality \eqref{keyi} guarantees that, for all $k\in \mathbb{N}$,
\begin{align*}
\left\| x_{k+1} - x^\star \right\|^2
%&\leq 
%\alpha_k \left\| x_{0} - x^\star \right\|^2 + \left\| x_k - x^\star \right\|^2 
%+ 2 \gamma_k \left( f^{(w_k)} (x^\star) - f^{(w_k)} (z_k) \right)\\
&\leq 
\alpha_k \left\| x_{0} - x^\star \right\|^2 + \left\| x_k - x^\star \right\|^2 
+ 2 \gamma_k \left( f^{(w_k)} (x^\star) - f^{(w_k)} (x_k) \right)\\
&\quad + 2 \gamma_k \left( f^{(w_k)} (x_k) - f^{(w_k)} (z_k) \right),
\end{align*}
which, together with the nonempty condition of $\partial f^{(w_k)}(x_k)$ and the triangle inequality, implies that, for all $k\in \mathbb{N}$, there exists $\bar{u}_k \in \partial f^{(w_k)}(x_k)$ such that
\begin{align*}
\left\| x_{k+1} - x^\star \right\|^2
&\leq 
\alpha_k \left\| x_{0} - x^\star \right\|^2 + \left\| x_k - x^\star \right\|^2 
+ 2 \gamma_k \left( f^{(w_k)} (x^\star) - f^{(w_k)} (x_k) \right)\\
&\quad + 2 \gamma_k \left\|\bar{u}_k \right\| \left\| x_k - z_k \right\|.
\end{align*}
Accordingly, 
\begin{align*}
\left\| x_{n+m+1} - x^\star \right\|^2
&\leq 
\left\| x_{n+1} - x^\star \right\|^2 
 + 2 \sum_{k=n+1}^{n+m} \gamma_k \left( f^{(w_k)} (x^\star) - f^{(w_k)} (x_k) \right)\\
&\quad + \left\| x_{0} - x^\star \right\|^2 \sum_{k=n+1}^{n+m} \alpha_k
+ 2 \sum_{k=n+1}^{n+m} \gamma_k \left\|\bar{u}_k \right\| \left\| x_k - z_k \right\|.
\end{align*}
A discussion similar to the one for obtaining \eqref{ineq:f} implies that
\begin{align}\label{Key}
\begin{split}
&\frac{2}{\gamma_{n+m}} \sum_{k=n+1}^{n+m} \gamma_k \left(f^{(w_k)} (x_k) - f^{(w_k)} (x^\star) \right)\\
\leq& \left(\left\| x_{n+1} - x^\star \right\| + \left\| x_{n+m+1} - x^\star \right\|\right) \frac{\left\| x_{n+m+1} - x_{n+1}\right\|}{\gamma_{n+m}}
+ \frac{\left\| x_{0} - x^\star \right\|^2}{\gamma_{n+m}} \sum_{k=n+1}^{n+m} \alpha_k\\
&+ \frac{2}{\gamma_{n+m}} \sum_{k=n+1}^{n+m} \gamma_k \left\|\bar{u}_k \right\| \left\| x_k - z_k \right\|.
\end{split}
\end{align}
Hence, the same discussion as for the proof of Lemma \ref{lem:4}, together with Lemma \ref{lem:3_1} and Assumptions \ref{omega} and \ref{stepsize2}, leads to the finding that $\limsup_{n\to\infty} f(x_n) \leq f^\star$ almost surely. Furthermore, the same discussion as for the proof of Theorem \ref{thm:1}, together with Assumption \ref{bounded2} and Lemmas \ref{lem:3_1} and \ref{lem:4_1}, guarantees that any weak sequential cluster point of $(x_n)_{n\in\mathbb{N}}$ almost surely belongs to $X^\star$. This means that Theorem \ref{thm:3} holds.

\subsection{Convergence rate analysis of Algorithm \ref{algo:2}} 
The following proposition establishes the rate of convergence for Algorithm \ref{algo:2}.

\begin{prop}\label{thm:4}
Suppose that the assumptions in Theorem \ref{thm:3} hold and that $(x_n)_{n\in \mathbb{N}}$ is the sequence generated by Algorithm \ref{algo:2}. Then, there exist $K_i \in \mathbb{R}$ ($i=1,2$) such that, for all $n\in \mathbb{N}$, almost surely
\begin{align*}
\left\| x_n - T^{(i)} (x_n) \right\|
\leq 
\sqrt{K_1 \alpha_{n} + K_2 \gamma_{n}}. 
\end{align*}
Moreover, under Sub-assumption \ref{omega}(iii), if there exists $k_0 \in \mathbb{N}$ such that almost surely $f (x_n) \geq f^\star$ for all $n \geq k_0$, there exist $k_1 \in \mathbb{N}$ and $K_i \in \mathbb{R}$ ($i=3,4,5,6$) such that, for all $n \geq \max\{k_0,k_1\}$, almost surely
\begin{align}\label{rate:2_2}
\frac{1}{m} \sum_{k=n+1}^{n+m} f(x_k) - f^\star
\leq
K_3 \frac{o(\gamma_{n+m})}{\gamma_{n+m}}
+ K_4\frac{\alpha_n}{\gamma_n}
+ \sqrt{K_5 \alpha_n + K_6 \gamma_n}.
\end{align}
Under Sub-assumption \ref{omega}(iv), there exist $K_i \in \mathbb{R}$ ($i=7,8,9,10$) such that, for all $n \in \mathbb{N}$, almost surely
\begin{align}\label{rate:2_3}
f(x_n) - f^\star 
&\leq 
K_7 \frac{o(\gamma_{n+m})}{\gamma_{n+m}}
+ K_8 \frac{\alpha_n}{\gamma_n} 
+ \sqrt{K_9 \alpha_n + K_{10} \gamma_n}.
\end{align}
\end{prop}

Let us compare Algorithm \ref{algo:1} with Algorithm \ref{algo:2}. Algorithm \ref{algo:1} can be applied to only smooth convex stochastic optimization whereas Algorithm \ref{algo:2} can be applied to nonsmooth convex stochastic optimization. Theorem \ref{thm:3} guarantees that any weak sequential cluster point of the sequence generated by Algorithm \ref{algo:2} almost surely belongs to the solution set of Problem \ref{prob:1} when $f^{(i)}$ $(i\in \mathcal{I})$ is smooth. Proposition \ref{thm:4} implies that Algorithm \ref{algo:2} with $\gamma_n := 1/n^a$ and $\alpha_n := 1/n^b$ $(n\geq 1)$, where $a \in (0,1/2)$, $b\in (a,1-a)$, and $a+b < 1$, satisfies, for all $i\in \mathcal{I}$,
\begin{align*}
\left\| x_n - T^{(i)} (x_n) \right\| 
= O \left(\frac{1}{\sqrt{n^{a}}} \right),
\end{align*}
which is the same as the rate of convergence of Algorithm \ref{algo:1} with $\lambda_n := 1/n^a$ and $\alpha_n := 1/n^b$ $(n\geq 1)$ for $\|x_n - T^{(i)}(x_n)\|$ (see \eqref{rate:1_1}). Moreover, \eqref{rate:2_2} and \eqref{rate:2_3} imply that, under the assumptions in Proposition \ref{thm:4} and the condition $o(\gamma_{n}) = 1/n^c$, where $c > a$, 
\begin{align*}
\frac{1}{m} \sum_{k=n+1}^{n+m} f(x_k) - f^\star
= 
O \left( \frac{1}{n^{\min\{ a/2,b-a,c-a \}}} \right), \text{ }
f(x_n) - f^\star 
= 
O \left( \frac{1}{n^{\min\{ a/2,b-a,c-a \}}} \right).
\end{align*}
Therefore, under the condition that $\lambda_n = \gamma_n := 1/n^a$ and $\alpha_n := 1/n^b$ $(n\geq 1)$, the rate of convergence of Algorithm \ref{algo:1} (see \eqref{Rate} and \eqref{Rate:1}) is almost the same as that of Algorithm \ref{algo:2}.

{\em Proof:}
Since \eqref{constant1} and \eqref{constant2} hold, the almost sure boundedness of $(x_n)_{n\in \mathbb{N}}$ and the monotone decreasing conditions of $(\alpha_n)_{n\in\mathbb{N}}$ and $(\gamma_n)_{n\in\mathbb{N}}$ with $\lim_{n\to\infty} \alpha_n = 0$ mean the existence of $K_i \in \mathbb{R}$ $(i=2,3)$ such that, for all $n\in \mathbb{N}$, almost surely
\begin{align*}
\mathbb{E}\left[ \left\| z_n - x_n \right\|^2 \Big| \mathcal{F}_n \right]
\leq K_2 \alpha_n + K_3 \gamma_n, \text{ } 
\mathbb{E}\left[ \left\| z_n - \mathsf{T}^{(w_n)}(z_n) \right\|^2 \Big| \mathcal{F}_n \right]
\leq K_2 \alpha_n + K_3 \gamma_n,
\end{align*}
which, together with \eqref{constant3} and the existence of $K_4\in \mathbb{R}$ such that, for all $i\in \mathcal{I}$ and for all $n\in \mathbb{N}$, almost surely $\| x_n - T^{(i)}(x_n)\|^2 \leq K_4 \mathbb{E}[\|x_n - \mathsf{T}^{(w_n)}(x_n)\|^2|\mathcal{F}_n]$ (see Sub-assumption \ref{omega}(ii)), means that, for all $i\in \mathcal{I}$ and for all $n\in \mathbb{N}$, almost surely
\begin{align*}
\left\| x_n - T^{(i)}(x_n)\right\|^2
\leq 10 K_4 \left(K_2 \alpha_n + K_3 \gamma_n \right).
\end{align*}
The same discussion as for obtaining \eqref{key} means that almost surely
\begin{align}\label{key_2}
\mathbb{E} \left[\frac{\left\| x_{n+m+1} - x_{n+1} \right\|}{\gamma_{n+m}} \Bigg| \mathcal{F}_{n} \right]
= \frac{o(\gamma_{n+m})}{\gamma_{n+m}}.
\end{align}

Let us assume that Sub-assumption \ref{omega}(iii) holds. Lemma \ref{lem:1_1}, the almost sure boundedness of $(x_n)_{n\in\mathbb{N}}$, and \eqref{constant1} imply that, for all $n\in \mathbb{N}$, almost surely $\| x_n - z_n \| \leq \sqrt{K_2 \alpha_n + K_3 \gamma_n}$. Taking the expectation in \eqref{Key} conditioned on $\mathcal{F}_n$ thus guarantees that there exist $K_i \in \mathbb{R}$ $(i=7,8,9)$ such that, for all $n \in \mathbb{N}$, almost surely
\begin{align*}
\frac{2}{\gamma_{n+m}} \sum_{k=n+1}^{n+m} \gamma_k \left( f(x_k) - f^\star \right)
&\leq
K_7 \frac{o(\gamma_{n+m})}{\gamma_{n+m}}
+ \frac{K_8}{\gamma_{n+m}}\sum_{k=n+1}^{n+m} \alpha_k
 + \frac{2 K_9}{\gamma_{n+m}} \sum_{k=n+1}^{n+m} \gamma_k
\sqrt{K_2 \alpha_k + K_3 \gamma_k},
\end{align*} 
where $K_9 := \max_{i\in \mathcal{I}} \sup \{ \| u_n^{(i)} \| \colon u_n^{(i)} \in \partial f^{(i)}(x_n), n\in \mathbb{N}\} < \infty$ comes from the almost sure boundedness of $(x_n)_{n\in \mathbb{N}}$ and Proposition \ref{prop:1}(iii). Accordingly, from the existence of $m_0 \in \mathbb{N}$ such that almost surely $f(x_n) - f^\star \geq 0$ for all $n \geq m_0$, \eqref{stepsizes}, and the monotone decreasing conditions of $(\alpha_n)_{n\in \mathbb{N}}$ and $(\gamma_n)_{n\in \mathbb{N}}$, there exist $K_{i} \in \mathbb{R}$ $(i=10,11)$ such that, for all $n \geq m_0$, almost surely
\begin{align*}
\frac{1}{m} \sum_{k=n+1}^{n+m} f(x_k) - f^\star
\leq
\frac{K_7}{m} \frac{o(\gamma_{n+m})}{\gamma_{n+m}}
+ \frac{K_{10}}{m}\frac{\alpha_n}{\gamma_n}
+ K_{11}\sqrt{K_2 \alpha_n + K_3 \gamma_n}.
\end{align*}

Next, let us assume that Sub-assumption \ref{omega}(iv) holds. From \eqref{Key} and the definition of $\partial f^{(w_k)}$, for all $n\in \mathbb{N}$, almost surely
\begin{align}\label{Key1}
\begin{split}
&2 \sum_{k=n+1}^{n+m} \left(f^{(w_k)} (x_n) - f^{(w_k)} (x^\star) \right)\\
\leq& \left(\left\| x_{n+1} - x^\star \right\| + \left\| x_{n+m+1} - x^\star \right\|\right) \frac{\left\| x_{n+m+1} - x_{n+1}\right\|}{\gamma_{n+m}}
+ \frac{\left\| x_{0} - x^\star \right\|^2 m \alpha_{n+m}}{\gamma_{n+m}}\\
&+ 2K_9 \sum_{k=n+1}^{n+m} \left\| x_k - z_k \right\| 
+ 2K_9 \sum_{k=n+1}^{n+m} \left\| x_n - x_k\right\|.
\end{split}
\end{align}
From \eqref{constant1} and \eqref{constant2}, the triangle inequality means that, for all $n\in \mathbb{N}$, almost surely
\begin{align*}
\left\| x_{n+1} - x_n\right\|
&\leq \left\| x_{n+1} - y_n\right\| + \left\| \mathsf{T}^{(w_n)}(z_n) - z_n \right\|
+ \left\| z_n - x_n \right\|\\
&\leq \alpha_n \left\|x_0 - y_n \right\| + 2 \sqrt{K_2 \alpha_n + K_3 \gamma_n},
\end{align*}
which, together with the triangle inequality and the monotone decreasing conditions of $(\alpha_n)_{n\in \mathbb{N}}$ and $(\gamma_n)_{n\in\mathbb{N}}$, means that, for all $n\in \mathbb{N}$, almost surely
\begin{align*}
\sum_{k=n+1}^{n+m} \left\| x_n - x_k\right\|
&\leq \sum_{j=0}^{m-1} \left(m-j\right) \left\| x_{n+j} - x_{n+j+1} \right\|\\
&\leq \sum_{j=0}^{m-1} \left(m-j\right) 
 \left( K_{11}\alpha_{n+j} + 2 \sqrt{K_2 \alpha_{n+j} + K_3 \gamma_{n+j}}\right)\\
&\leq \frac{m(m+1)}{2}
 \left( K_{11}\alpha_{n} + 2 \sqrt{K_2 \alpha_{n} + K_3 \gamma_{n}}\right),
\end{align*}
where almost surely $K_{11} := \sup \{ \|x_0 - y_n\| \colon n\in \mathbb{N} \} < \infty$. Taking the expectation in \eqref{Key1} conditioned on $\mathcal{F}_n$ thus guarantees that, for all $n \in \mathbb{N}$, almost surely
\begin{align*}
m\left(f(x_n) - f^\star \right)
&\leq 
K_7 \frac{o(\gamma_{n+m})}{\gamma_{n+m}}
+ K_{10} \frac{\alpha_n}{\gamma_n}\\
&\quad + m K_8 \sqrt{K_2 \alpha_n + K_3 \gamma_n}
+ \frac{m(m+1)}{2} K_8 \left( K_{11}\alpha_{n} + 2 \sqrt{K_2 \alpha_{n} + K_3 \gamma_{n}}\right),
\end{align*}
which, together with (C4), completes the proof.

\begin{rem}
{\em A discussion similar to the one for obtaining \eqref{key_1} ensures that there exist $m_0 \in \mathbb{N}$ and $\bar{K}_i \in \mathbb{R}$ $(i=1,2)$ such that, for all $n > m_0$, almost surely
\begin{align}\label{key_3}
\mathbb{E} \left[\frac{\left\| x_{n+m(n)+1} - x_{n+1} \right\|}{\gamma_{n+m(n)}} \Bigg| \mathcal{F}_{n} \right]
\leq
\bar{K}_1 \prod_{k=m_0}^{n} \left(1 - \alpha_{k+m(k)} \right) + \bar{K}_2 N(n),
\end{align}
where $N(n):= \max \{ (1 /\gamma_{k+m(k)}) | 1 - \alpha_k/\alpha_{k+m(k)}|, (1/\alpha_{k+m(k)}) | (1/\gamma_{k+m(k)}) - (1/\gamma_{k+m(k)-1}) |, (1/\alpha_{k+m(k)}) |\gamma_{k+m(k)} - \gamma_k|/\gamma_{k+m(k)}^2 \colon k=n,n-1,\ldots,m_0 \}$. Accordingly, \eqref{key_2} can be replaced with \eqref{key_3}.
}
\end{rem}

\section{Numerical results}
\label{sec:5}
This section considers Problem \ref{prob:1} when $f^{(i)} \colon \mathbb{R}^d \to \mathbb{R}$ and $T^{(i)} \colon \mathbb{R}^d \to \mathbb{R}^d$ $(i\in \mathcal{I})$ are defined for all $x := (x_1, x_2,\ldots,x_d) \in \mathbb{R}^d$ by
\begin{align}
&f^{(i)}(x) := 
\frac{1}{2} \left\langle x, A^{(i)}x \right\rangle + \left\langle b^{(i)},x\right\rangle
\text{ or }
\sum_{j\in \mathcal{D}}\omega_{j}^{(i)}\left| x_j - a_j^{(i)}\right|,\nonumber\\
&T^{(i)}(x) := \frac{1}{2} \left[ x 
 + P_C \left(\frac{1}{K} \sum_{k\in \mathcal{K}} P_{C_k^{(i)}} (x)\right)
 \right],\label{concrete}
\end{align}
where $A^{(i)} \in \mathbb{R}^{d \times d}$ is a diagonal matrix with diagonal components $\lambda_j^{(i)} \geq 0$, $b^{(i)} \in \mathbb{R}^d$, $\omega_j^{(i)} > 0$, $a_j^{(i)} \in \mathbb{R}$, $r_k^{(i)} > 0$, $c_k^{(i)} \in \mathbb{R}^d$, $C_k^{(i)} := \{ x\in \mathbb{R}^d \colon \| x - c_k^{(i)} \| \leq r_k^{(i)} \}$ $(i\in \mathcal{I}, k\in \mathcal{K}:= \{1,2,\ldots,K\},j\in \mathcal{D} := \{1,2,\ldots,d\})$, and $C := \{ x\in \mathbb{R}^d \colon \|x\| \leq 1\}$.

Since the metric projection onto each of $C$ and $C_k^{(i)}$ $(i\in \mathcal{I},k\in \mathcal{K})$ can be computed within a finite number of arithmetic operations, $T^{(i)}$ $(i\in \mathcal{I})$ defined by \eqref{concrete} can be computed efficiently. Moreover, $T^{(i)}$ $(i\in \mathcal{I})$ satisfies the firm nonexpansivity condition (see (A1)), and $\mathrm{Fix}(T^{(i)})$ coincides with a subset of $C$ with the elements closest to $C_k^{(i)}$s in terms of the mean square norm \cite[Proposition 4.2]{yamada}. This subset, denoted by $C_{\Phi}^{(i)} := \{x\in C \colon \Phi^{(i)}(x) : = (1/K) \sum_{k\in \mathcal{K}} (\min_{z\in C_k^{(i)}} \|x - z\|)^2 = \min_{y\in C} \Phi^{(i)}(y)\}$ $(= \mathrm{Fix}(T^{(i)}))$, is called the {\em generalized convex feasible set} \cite{com1999,yamada}, which is well defined even when $C \cap \bigcap_{k\in \mathcal{K}} C_k^{(i)} = \emptyset$ (see \cite{com1999,iiduka_siam2012,iiduka_yamada_siam2012,yamada} for applications of the generalized convex feasible set). The boundedness of $C$ guarantees that $\mathrm{Fix}(T^{(i)}) = C_{\Phi}^{(i)} \neq \emptyset$ \cite[Remark 4.3(a)]{yamada}. 

The experimental evaluations of the two proposed algorithms were done using a Mac Pro with a 3-GHz 8-Core Intel Xeon E5 processor and 32-GB 1866-MHz DDR3 memory. The algorithms were written in Java (version 9) with $d := 2^{10} = 1024$, $I := 16$, and $K := 3$. The values of $\lambda_j^{(i)} \in [0,d]$, $b^{(i)} \in [-1,1]^d$, $\omega_j^{(i)} \in (0,1]$, $a_j^{(i)} \in [-1,1]$, $r_k^{(i)} \in (0,1]$, and $c_k^{(i)} \in [-1/\sqrt{d}, 1/\sqrt{d})^d$ were randomly generated using the Mersenne Twister pseudorandom number generator (provided by Apache Commons Math 3.6). Algorithm \ref{algo:1} (resp. Algorithm \ref{algo:2}) was used with \eqref{y_n} (resp. \eqref{y_n2}), which implies the boundedness of $(y_n)_{n\in\mathbb{N}}$ (see Assumptions \ref{bounded} and \ref{bounded2}). The step-size sequences were $\lambda_n = \gamma_n := 10^{-3}/(n+1)^a$ and $\alpha_n := 10^{-3}/(n+1)^b$, where $(a,b)$ is (A) $(1/4,1/2)$ or (B) $(1/8,3/4)$, which satisfy Assumptions \ref{stepsize} and \ref{stepsize2}.\footnote{Existing fixed point optimization algorithms \cite{iiduka2016,iiduka_ejor2016} with small step sizes (e.g., $\gamma_n := 10^{-2}/(n+ 1)^a, 10^{-3}/(n+1)^a$) have faster convergence. Hence, the experiment used step sizes $\lambda_n = \gamma_n := 10^{-3}/(n+1)^a$ and $\alpha_n := 10^{-3}/(n+1)^b$.}

To see how the choice of $(w_n)_{n\in\mathbb{N}}$ affects the convergence rate of the two algorithms, Algorithms \ref{algo:1} and \ref{algo:2} were used with one of the following conditions.
\begin{enumerate}
\item[(I)] The samples were generated nearly independently; i.e., for all $i\in \mathcal{I}$, there existed $\rho_i \in (0,1]$ such that, almost surely $\inf_{n\in \mathbb{N}} \mathbb{P}(w_n = i | \mathcal{F}_n ) \geq \rho_i/I$.
\item[(II)] The samples were selected to be nonexpansive mappings of which the fixed point sets were the most distant from the current iterates; i.e., $w_n \in \argmax_{i\in \mathcal{I}} \| x_n - T^{(i)}(x_n) \|^2$ for all $n\in \mathbb{N}$.
\item[(III)] The samples were generated in accordance with a random permutation of the indexes within a cycle; i.e., for all $t \in \mathbb{N}$, $(\mathsf{T}^{(w_n)})_{n\in\mathbb{N}}$, where $n=tI, tI+1, \ldots, (t+1)I-1$, was a permutation of $\{ T^{(1)}, T^{(2)}, \ldots, T^{(I)} \}$.
\item[(IV)] The samples were generated through state transitions of a Markov chain; i.e., $(w_n)_{n\in \mathbb{N}}$ was generated using an irreducible and aperiodic Markov chain with states $1,2,\ldots,I$.
\end{enumerate} 

Conditions (I)--(IV) were defined on the basis of Assumptions 4--7 in \cite{bert2015}. The conclusions in \cite{bert2015} show that the sequence $(w_n)_{n\in\mathbb{N}}$ in each condition satisfies Sub-assumptions \ref{omega}(i) and (ii). In the experiment, $(w_n)_{n\in\mathbb{N}}$ in (IV) was generated using a positive Markov matrix with randomly chosen elements.

One hundred samplings, each starting from a different randomly chosen initial point, were performed, and the results were averaged. Two performance measures were used. For each $n\in\mathbb{N}$,
\begin{align*}
D_n := \frac{1}{100} \sum_{s=1}^{100} \sum_{i\in\mathcal{I}} \left\| x_n(s) - T^{(i)} \left(x_n(s) \right) \right\|
\text{ and }
F_n := \frac{1}{100} \sum_{s=1}^{100}
\mathbb{E} \left[f^{(w)} \left(x_n (s) \right) \right],
\end{align*}
where $(x_n(s))_{n\in\mathbb{N}}$ is the sequence generated from initial point $x(s)$ $(s=1,2,\ldots,100)$ for each of the two algorithms. The value of $D_n$ represents the mean value of the sums of the distances between $x_n(s)$ and $T^{(i)}(x_n(s))$. Hence, if $(D_n)_{n\in\mathbb{N}}$ converges to $0$, $(x_n)_{n\in\mathbb{N}}$ converges to some point in $\bigcap_{i\in\mathcal{I}} \mathrm{Fix}(T^{(i)}) = \bigcap_{i\in\mathcal{I}} C_{\Phi}^{(i)}$.  
$F_n$ is the average of $\mathbb{E} [f^{(w)}(x_n (s))]$ $(s=1,2,\ldots,100)$, and the values of $F_n$ generated by Algorithms \ref{algo:1} and \ref{algo:2} with Conditions (I)--(IV) differ since the samples are coming from different distributions in (I)--(IV).  
The stopping condition was $n=1000$.
 
First, let us consider the problem when $f^{(i)}(x):= (1/2) \langle x,A^{(i)}x \rangle + \langle b^{(i)},x\rangle$ $(i\in \mathcal{I})$ (i.e., $f^{(i)}$ is smooth and convex), which can be solved using Algorithm \ref{algo:1}. Table \ref{table1} shows the number of iterations $n$ and elapsed time when Algorithm \ref{algo:1} with one of (I)--(IV) and one of (A) and (B) satisfied $D_n \leq 10^{-3}$ and $|F_n - F_{n-1}| \leq 10^{-5}$. All the algorithms converged to a point in $\bigcap_{i\in \mathcal{I}} \mathrm{Fix}(T^{(i)})$ in the early stages. $F_n$ when Algorithm \ref{algo:1} satisfied $|F_n - F_{n-1}| \leq 10^{-5}$ was different from $F_{1000}$ because the behavior of Algorithm \ref{algo:1} was unstable in the early stages. Checking showed that Algorithm \ref{algo:1} satisfied $D_{n} \approx 0$ for $n \geq 10$ and that its behavior was stable for $n \geq 900$. When one of (I)--(IV) was fixed, $F_{1000}$ generated by Algorithm \ref{algo:1}(A) was smaller than $F_{1000}$ generated by Algorithm \ref{algo:1}(B). Accordingly, Algorithm \ref{algo:1}(A) performed better than Algorithm \ref{algo:1}(B). 

{\small
\begin{table}
% \centering
 \caption{Behavior of $D_n$ and $F_n$ for Algorithm \ref{algo:1}}
 \label{table1}
 \begin{tabular}{l||ccc|ccc|cc}
 \hline
 & \multicolumn{3}{c|}{$D_n \leq 10^{-3}$} 
 & \multicolumn{3}{c|}{$|F_{n}-F_{n-1}| \leq 10^{-5}$}
 & \multicolumn{2}{c}{$n=1000$}\\
 \hline
 & $n$ & time [s] & $D_n$ & $n$ & time [s] & $F_n$ & time [s] & $F_n$\\
 \hline
Alg.\ref{algo:1}(I)(A) 
&6 & 0.000071 & 0.000392 
& 132 & 0.001263 & 0.022259 
&0.009604 & $-0.011421$\\
Alg.\ref{algo:1}(I)(B) 
&6 & 0.000072 & 0.000000 
& 301 & 0.002930 & 0.010067 
& 0.009704 & 0.002931\\
Alg.\ref{algo:1}(II)(A) 
&6 & 0.001134 & 0.000000 
& 250 & 0.037588 & 0.024567 
& 0.147227 & 0.006503\\
Alg.\ref{algo:1}(II)(B) 
&5 & 0.000963 & 0.000100 
& 99 & 0.015022 & 0.038290
& 0.147379 & 0.019506\\
Alg.\ref{algo:1}(III)(A) 
&5 & 0.000061 & 0.000000 
& 78 & 0.000754 & 0.033988 
& 0.009607 & $-0.005376$\\
Alg.\ref{algo:1}(III)(B) 
&4 & 0.000051 & 0.000000 
& 110 & 0.001063 & 0.019045 
& 0.009731 & 0.003758\\
Alg.\ref{algo:1}(IV)(A) 
&5 & 0.000062 & 0.000065 
& 423 & 0.004195 & 0.016480
& 0.009871 & 0.007351\\
Alg.\ref{algo:1}(IV)(B) 
&5 & 0.000062 & 0.000000
& 484 & 0.004830 & 0.025469
& 0.009889 & 0.020530 \\
 \hline
 \end{tabular}
\end{table}
}

Next, let us consider the case in which $f^{(i)}(x):= \sum_{j\in \mathcal{D}} \omega_j^{(i)}|x_j - a_j^{(i)}|$ $(i\in \mathcal{I})$ (i.e., $f^{(i)}$ is nonsmooth and convex), which can be solved using Algorithm \ref{algo:2}. Table \ref{table2} shows that all the algorithms optimized $F_n$ in the early stages and then searched for a point in $\bigcap_{i\in \mathcal{I}} \mathrm{Fix}(T^{(i)})$, in contrast to Algorithm \ref{algo:1} (see Table \ref{table1}). Checking showed that the behavior of Algorithm \ref{algo:2} was stable. Moreover, when one of (I)--(IV) was fixed, Algorithm \ref{algo:2}(B) satisfied $D_n \leq 10^{-2}$ more quickly than Algorithm \ref{algo:2}(A), and $F_{1000}$ generated by Algorithm \ref{algo:2}(B) was smaller than that generated by Algorithm \ref{algo:2}(A). Accordingly, Algorithm \ref{algo:2}(B) performed better than Algorithm \ref{algo:2}(A).

{\small
\begin{table}
% \centering
 \caption{Behavior of $D_n$ and $F_n$ for Algorithm \ref{algo:2}}
 \label{table2}
 \begin{tabular}{l||ccc|ccc|cc}
 \hline
 & \multicolumn{3}{c|}{$D_n \leq 10^{-2}$} 
 & \multicolumn{3}{c|}{$|F_{n}-F_{n-1}| \leq 10^{-5}$}
 & \multicolumn{2}{c}{$n=1000$}\\
 \hline
 & $n$ & time [s] & $D_n$ & $n$ & time [s] & $F_n$ & time [s] & $F_n$\\
 \hline
Alg.\ref{algo:2}(I)(A) 
&$>1000$ & --- & --- 
& 14 & 0.000191 & 0.206366
&0.010437 & 0.202815\\
Alg.\ref{algo:2}(I)(B) 
&522 & 0.005282 & 0.009996 
& 14 & 0.000193 & 0.206323
&0.009920 & 0.194252\\
Alg.\ref{algo:2}(II)(A) 
&770 & 0.120310 & 0.009993 
& 9 & 0.001693 & 0.193260
&0.155903 & 0.191289\\
Alg.\ref{algo:2}(II)(B) 
&46 & 0.007322 & 0.009871
& 9 & 0.001651 & 0.193248
& 0.148430 & 0.187172\\
Alg.\ref{algo:2}(III)(A) 
&771 & 0.008040 & 0.009961 
&14 & 0.000194 & 0.193133
&0.010388 & 0.191453\\
Alg.\ref{algo:2}(III)(B) 
&96 & 0.001040 & 0.009769
& 14 & 0.000190 & 0.193123 
& 0.009999 & 0.187654\\
Alg.\ref{algo:2}(IV)(A) 
&976 & 0.010334 & 0.009998 
&7 & 0.000106 & 0.193286
&0.010596 & 0.191582\\
Alg.\ref{algo:2}(IV)(B) 
&121 & 0.001304 & 0.009790 
& 7 & 0.000109 & 0.193281
& 0.009994 & 0.188013 \\
 \hline
 \end{tabular}
\end{table}
}

\section{Conclusion}
\label{sec:6}
Two stochastic optimization algorithms were proposed for solving the problem of minimizing the expected value of convex functions over the intersection of fixed point sets of nonexpansive mappings in a real Hilbert space. One algorithm blends a stochastic gradient method with the Halpern fixed point algorithm while the other is based on a stochastic proximal point algorithm and the Halpern fixed point algorithm. Consideration of a case in which the step-size sequences are diminishing demonstrated that any weak sequential cluster point of the sequence generated by each of the two algorithms almost surely belongs to the solution set of the problem under certain assumptions. Convergence rate analysis of the two algorithms illustrated their efficiency. A discussion of concrete convex optimization over fixed point sets and the numerical results demonstrated their effectiveness.

{\small \textbf{Acknowledgments.} 
I am sincerely grateful to the Senior Editor, Takashi Tsuchiya, the anonymous associate editor, and the anonymous reviewers for their insightful comments. 
I also thank Kazuhiro Hishinuma for his input on the numerical examples. }

\end{document}